%% file: manuscript-jnsao.tex
\documentclass{jnsao}

\usepackage{thmtools}
\usepackage{zref-clever}
\usepackage{zref-vario}
\usepackage{algpseudocode}

\input{macros.sty}

\title{The Riemannian Convex Bundle Method}

\author{%
Ronny Bergmann%
\thanks{Norwegian University of Science and Technology, Department of Mathematical Sciences, NO-7041 Trondheim, Norway, \email{\detokenize{ronny.bergmann@ntnu.no}}, \url{https://www.ntnu.edu/employees/ronny.bergmann}, \orcid{0000-0001-8342-7218}}
\and
Roland Herzog%
\thanks{Interdisciplinary Center for Scientific Computing, Heidelberg University, 69120 Heidelberg, Germany \emph{and}
Institute for Mathematics Heidelberg University, 69120 Heidelberg, Germany,
\email{\detokenize{roland.herzog@iwr.uni-heidelberg.de}}, \url{https://scoop.iwr.uni-heidelberg.de}, \orcid{0000-0003-2164-6575}}%
\and
Hajg Jasa%
\thanks{Simula Research Laboratory, Department of Numerical Analysis and Scientific Computing, 0164 Oslo, Norway, \email{\detokenize{hajg@simula.no}}, \url{https://www.simula.no/people/hajg}, \orcid{0009-0002-7917-0530}}
}
\shortauthor{Bergmann, Herzog, Jasa}

\manuscriptsubmitted{2025-07-25}
\manuscriptaccepted{2026-09-22}
\manuscriptvolume{6}
\manuscriptnumber{16115}
\manuscriptyear{2026}
\manuscriptdoi{10.46298/jnsao-2026-16115}

\begin{document}
\maketitle

\begin{abstract}
\input{abstract.tex}
\end{abstract}

\input{main.tex}

\appendix
\input{appendix.tex}

\bibliographystyle{jnsao}
\bibliography{manuscript-jnsao}

\end{document}

%% file: abstract.tex
We introduce the convex bundle method to solve convex, non-smooth optimization problems on Riemannian manifolds of bounded sectional curvature.
Each step of our method is based on a model that involves the convex hull of previously collected subgradients, parallelly transported to the tangent space at the current serious iterate.
This generalizes the dual approach to nonsmooth optimization via classical bundle subproblems in Euclidean space.
We prove that, under mild conditions, the convex bundle method converges to a minimizer.
Several numerical examples implemented using \manoptjl illustrate the performance of the proposed method and compare it to the subgradient method, the cyclic proximal point algorithm, as well as the proximal bundle method.
We showcase how the proposed method outperforms these algorithms in terms of CPU runtime and number of iterations on examples for finding Riemannian medians, while it is comparable in performance on large signal denoising tasks.

%% file: main.tex
\section{Introduction}
\label{section:introduction}

Research on optimization problems defined on Riemannian manifolds has experienced steadily increasing activity in the past decades.
Among these, non-smooth, unconstrained optimization defines a class of particularly interesting problems.
For instance, tasks such as image denoising~\cite{BacakBergmannSteidlWeinmann:2016:1,LellmannStrekalovskiyKoetterCremers:2013:1,WeinmannDemaretStorath:2014:1}, inpainting~\cite{BergmannFitschenPerschSteidl:2018:1}, and interpolation~\cite{BergmannGousenbourger:2018:2,Zimmermann:2020:1} are often expressed in terms of high-dimensional, non-smooth optimization tasks defined on Riemannian manifolds.

We consider problems of the form
\begin{equation}
	\label{problem:main}
  \zcsetup{reftype=problem}
	\text{Minimize }
	f(p)
	,
	\quad
	p \in \cM
	,
\end{equation}
where $\cM$ is a Riemannian manifold and $f$ is geodesically convex and potentially non-smooth.
Our setting will be that of a manifold~$\cM$ with bounded sectional curvature.
The precise requirements are given in the beginning of \zcref{section:approximation-of-the-epsilon-subdifferential}.

In the Euclidean case, various algorithms lend themselves to the minimization of non-smooth convex functions with additional structure.
We mention the Douglas-Rachford splitting algorithm~\cite{DouglasRachford:1956:1}, the Alternating Direction Method of Multipliers (ADMM)~\cite{GabayMercier:1976:1}, and primal-dual splitting methods such as the Chambolle-Pock algorithm~\cite{ChambollePock:2011:1}.
On Riemannian manifolds, some of the first non-smooth optimization algorithms include the subgradient method~\cite{FerreiraOliveira:1998:1} and the proximal point algorithm~\cite{FerreiraOliveira:2002:1}.
Further algorithmic options exist on Hadamard manifolds, \eg, adaptations of Douglas-Rachford splitting~\cite{BergmannPerschSteidl:2016:1} or, quite recently, the Chambolle-Pock algorithm \cite{BergmannHerzogSilvaLouzeiroTenbrinckVidalNunez:2021:1} based on a newly introduced concept of duality.

A class of algorithms that have proved to be effective in optimizing generic, convex and non-convex, non-smooth functions in the Euclidean setting is that of bundle methods; see, \eg, \cite{Maekelae:2002:1} for a survey of bundle methods, and \cite{BonnansGilbertLemarechalSagastizabal:2006:1,GeigerKanzow:2002:1,Noll:2013:1} for textbook presentations on variants thereof.
These algorithms are characterized by the fact that the descent strategy is employed within a framework that guarantees stability by keeping track of past iterates and subgradients.
One subcategory of stable bundle methods is that of proximal bundle algorithms, where a proximal term is added to the model of the objective function.
A Riemannian bundle method for possibly non-convex, Lipschitz functions, that makes use of a proximal term and models the objective function on the tangent space was presented in~\cite{HoseiniMonjeziNobakhtianPouryayevali:2023:1}.
To obtain a search direction, this method solves a proximal subproblem that is formally different from the subproblem we solve in the convex case, most prominently due to the introduction of a proximal parameter.
The next candidate point is then obtained through the search direction.
This search direction may have to be adapted according to a restricted-step procedure that depends on the injectivity radius of the manifold.
A related method is the Riemannian restricted memory quasi-Newton bundle method \cite{TangXingHuangJian:2026:1} for minimizing locally Lipschitz functions.
This method approximates the curvature information of the objective function by applying a Riemannian quasi-Newton update.
A subgradient aggregation technique is then used to reduce the computations in the quadratic subproblem for procuring the candidate descent direction.
The stepsize selection entails a line-search that leverages a Riemannian semismoothness concept.

We are further interested in numerical applications and examples for non-smooth optimization problems on Riemannian manifolds, such as typical signal and manifold-valued image denoising problems.
In this context, in parallel to the aforementioned algorithmic advancements, ample software has been developed.
Examples include the \manopt family, consisting of the inaugural \matlab version \manopt~\cite{BoumalMishraAbsilSepulchre:2014:1}, the \python version \pymanopt~\cite{TownsendKoepWeichwald:2016:2}, and the \julia version \manoptjl~\cite{Bergmann:2022:1}.
Another popular software related to differential geometry and to some extent to optimization is the \python package \geomstats~\cite{MiolaneGuiguiLeBrigantMatheHouThanwerdasHeyderPeltreKoepZaatitiHajriCabanesGeraldChauchatShewmakeBrooksKainzDonnatHolmesPennec:2020:1}.
In this work we focus on the \julia package \manoptjl.

\subsection*{Contributions}

The aim of this paper is to generalize the \emph{convex} bundle method (CBM) for convex, non-smooth functions from Euclidean spaces to Riemannian manifolds with bounded sectional curvature.
A Euclidean version of this method is described, \eg, in~\cite[Chapter~6.7]{GeigerKanzow:2002:1}.
We generalize this method by exploring properties of geodesically convex functions, for which a Riemannian generalization of the convex subdifferential (also known as the Fenchel subdifferential) is defined.
At its core, our method approximates the $\varepsilon$-subdifferential of the objective function at the current serious iterate by the convex hull of the subgradients that were collected at the previous iterates.
This is shown in~\zcref{theorem:eps-approx}, where the main challenge arises from the fact that these subgradients lie in different tangent spaces and need to be parallely transported, which impacts their ability to be part of an approximate subdifferential at the target point.
The search direction in our algorithm is obtained by solving a dual subproblem arising from an orthogonal projection of the zero tangent vector onto the aforementioned inner approximation of the convex $\varepsilon$-subdifferential.
To procure the next candidate point, the search direction is shortened if necessary.
This is done in accordance with the function's effective domain, using a backtracking procedure.
If the manifold is strictly positively curved, a restriction of the domain of the objective may be needed for this set to be geodesically convex.
On the other hand, if the manifold is strictly negatively curved, we require the objective's domain to be bounded to prove that the sequence of serious iterates is bounded in~\zcref{lemma:bounded-min-sequence}.
\zcref{theorem:cbm-convergence} provides the same type of result as its Euclidean counterpart presented, \eg, in~\cite[Theorem~6.80]{GeigerKanzow:2002:1}, under the additional assumption that the domain of the objective be bounded whenever the manifold is not flat.
The convergence to a (global) minimizer of the objective is guaranteed whenever the set of interior solutions to \zcref{problem:main} is nonempty, and the sequence of serious iterates produced by the proposed algorithm does not have accumulation points on the boundary of $\dom f$.
This result hinges on the following observation.
Inspired by \cite[p.15]{BagirovGaudiosoKarmitsaMakelaTaheri:2020:1} and motivated in \zcref{remark:null-step-condition}, in case the algorithm produces a null step, we look for a step size that satisfies a geometric condition involving the linearization errors and curvature-dependent remainder terms defined in \zcref{eq:new-lin-errors} and \zcref{eq:fixed-remainder}, respectively.
This ensures that the last serious iterate output by our method is a minimizer in case the algorithm produces only a finite number of serious steps.

In summary, we
\begin{itemize}[nosep]
  \item approximate the $\varepsilon$-subdifferential in \zcref{theorem:eps-approx};
  \item introduce \zcref{algorithm:RCBM} and analyze its basic properties;
  \item show the asymptotic convergence of this method in the case where infinitely many serious steps occur in \zcref{lemma:infinite-opt-convergence};
  \item show the asymptotic convergence the proposed algorithm in the case where finitely many serious iterates are produced in \zcref{lemma:finite-opt-convergence};
  \item prove the full convergence of \zcref{algorithm:RCBM} in \zcref{theorem:cbm-convergence}.
\end{itemize}

\subsection*{Organization}

The remainder of this paper is organized as follows.
In \zcref{section:preliminaries} we recall some notions from Riemannian geometry and fix the notation.
\zcref{section:approximation-of-the-epsilon-subdifferential} introduces the main problem, the convex bundle method model, and an approximation to the $\varepsilon$-subdifferential tailored to Riemannian manifolds with bounded sectional curvature.
The derivation and convergence analysis of our method is motivated and inspired by~\cite[Chapter~6.7]{GeigerKanzow:2002:1}.
Proofs that directly carry over to the Riemannian case without any sort of modification are omitted, whereas proofs that require additional steps will be given.
\zcref{section:algorithm} is dedicated to the description of the algorithm.
\zcref{section:convergence} shows the analysis of the convergence properties of the proposed algorithm.
The same observations regarding the Euclidean case made for \zcref{section:approximation-of-the-epsilon-subdifferential} apply to this section as well.
In \zcref{section:numerics} we illustrate the performance of the Riemannian convex bundle method (RCBM) on several numerical examples on Riemannian manifolds.
\zcref{section:conclusion} concludes the paper.

\section{Preliminaries}
\label{section:preliminaries}

In this section, we recall some notions from Riemannian geometry.
For more details, the reader may wish to consult \cite{Lee:1997:1,DoCarmo:1992:1}.
Let $\cM$ be a smooth manifold with Riemannian metric $\riemannian{\cdot}{\cdot}$ and let $\covariantDerivativeSymbol$ be its Levi-Civita connection.
Given $p, q \in \cM$, we denote by
\begin{equation*}
	C_{pq}
	=
	\setDef[auto]{\gamma \colon[0,1] \to \cM}{\gamma \text{ is piecewise smooth and } \gamma(0) = p, \; \gamma(1) = q}
\end{equation*}
the set of all piecewise smooth curves from~$p$ to~$q$.
The \emph{arc length} of $\gamma \in C_{pq}$, $L \colon C_{pq} \to \bbR_{\ge 0}$ is defined as
\begin{equation*}
	L(\gamma)
	=
	\int_0^1 \riemanniannorm{\dot \gamma}[\gamma(t)] \, \d t
	,
\end{equation*}
where $\riemanniannorm{\cdot}[p]$ is the norm induced by the Riemannian metric in the tangent space $\tangentSpace{p}$ at the point~$p$.
The inner product given by the Riemannian metric at~$p$ is $\riemannian{\cdot}{\cdot}[p] \colon \tangentSpace{p} \times \tangentSpace{p} \to \bbR$.
The explicit reference to the base point will be omitted whenever the base point is clear from context.
Tangent vectors in $\tangentSpace{p}$ will be generically denoted by $X_p$ or $Y_p$ \etc where the subindex merely serves to indicate the base point.

The sectional curvature at $p \in \cM$ for linearly independent tangent vectors $X_p, Y_p \in \tangentSpace{p}$ is
\begin{equation}
	K_p(X_p,Y_p)
	=
	\frac{
		\riemannian{\mathrm{R}_p(X_p,Y_p) Y_p}{X_p}
	}{
		\riemanniannorm{X_p}^2
		\,
		\riemanniannorm{Y_p}^2
		-
		\riemannian{X_p}{Y_p}^2
	}
	,
\end{equation}
where $\mathrm{R}$ is the Riemann curvature tensor.

The Riemannian distance between two points, $\dist \colon \cM \times \cM \to \bbR_{\ge 0} \cup \{+\infty\}$, is given by
\begin{equation*}
	\dist(p,q)
	=
	\inf_{\gamma \in C_{pq}} L(\gamma)
	.
\end{equation*}
A \emph{geodesic arc} from $p \in \cM$ to $q \in \cM$ is a smooth curve $\gamma \in C_{pq}$ that is parallel along itself, \ie, $\covariantDerivative{\dot \gamma(t)}[\dot \gamma(t)] = 0$ for all $t \in [0,1]$.
A finite-length geodesic arc is \emph{minimal} if its arc length coincides with the Riemannian distance between its extreme points.
Given two points $p, q \in \cM$, we denote by $\geodesic{p}{q}$ a minimal geodesic arc that connects $p = \geodesic{p}{q}(0)$ to $q = \geodesic{p}{q}(1)$.

On Riemannian manifolds, there exist several notions of convex sets, which differ in strength; see for instance \cite[Section~2]{Walter:1974:1} or \cite[Definition~2.1]{WangLiYao:2016:1}.
We will only speak of sets referred to as \emph{strongly} convex sets in the literature.

\begin{definition}
	A set $\cU \subseteq \cM$ is said to be \emph{strongly convex} if for any $p,q \in \cU$, there exists a unique minimal geodesic arc $\geodesic{p}{q}$, and, moreover, $\geodesic{p}{q}$ lies entirely in~$\cU$.
\end{definition}

The exponential map at $p \in \cM$ is denoted by $\exponential{p} \colon \cU_p \to \cM$, where $\cU_p \subseteq \tangentSpace{p}$ is a suitable neighborhood of $0_p \in \tangentSpace{p}$
that depends on the injectivity radius of $\cM$ at the point $p$, which is the quantity given by
\begin{equation*}
	i_\cM(p)
	\coloneqq
	\sup
	\setDef[auto]{
		r > 0
	}{
		\exponential{p} \colon \openBall{r}{0_p} \to \exponential{p}(\openBall{r}{0_p})
		\text{ is a diffeomorphism}
	}%
	,
\end{equation*}
where $\openBall{r}{0_p} \subseteq \tangentSpace{p}$; see, \eg, \cite[Chapter~6, p.165]{Lee:2018:1}.
The (global) injectivity radius of $\cM$ is $i_\cM \coloneqq \inf_{p \in \cM} i_\cM(p)$.
The inverse of the exponential map at $p\in\cM$ is the logarithmic map, and it is denoted by $\logarithm{p}$.
It might not be defined on all of $\cM$, but it is guaranteed to be defined on some open ball of radius~$r \le i_\cM(p)$ centered at~$p$, denoted by $\openBall{r}{p} = \setDef[auto]{q \in \cM}{\dist(p,q) < r}$.
Note that given a strongly convex set~$\cU$ and $p, q \in \cU$, $\logarithm{p} q$ is the unique tangent vector that realizes the initial velocity of the unique minimal geodesic arc $\geodesic{p}{q}$ that lies in $\cU$.

We denote the \emph{parallel transport} from~$p$ to~$q$ along~$\geodesic{p}{q}$ with respect to the connection~$\covariantDerivativeSymbol$ by $\parallelTransport{p}{q} \colon \tangentSpace{p} \to \tangentSpace{q}$.
It is defined by
\begin{equation*}
	\parallelTransport{p}{q} Y_p
	=
	X_q
	\quad
	\text{for }
	Y_p \in \tangentSpace{p}
	,
\end{equation*}
where $X$ is the unique smooth vector field along $\geodesic{p}{q}$ which satisfies
\begin{equation*}
	X_p
	=
	Y_p
	\text{ and }
	\covariantDerivative{\dot \gamma}{X}
	=
	0
	.
\end{equation*}
Notice that $\parallelTransport{p}{q}$ depends on the choice of a minimal geodesic $\geodesic{p}{q}$ but in our setting it will not cause any ambiguity since $\geodesic{p}{q}$ is unique in a strongly convex set.
Denoting by $\parallelTransport{q}{p}$ the parallel transport along the curve $\geodesic{p}{q}$ traveled in the opposite direction, we have
\begin{equation}
	\label{eq:adjoint-of-parallel-transport}
	\riemannian{\parallelTransport{p}{q} Y_p}{Z_q}
	=
	\riemannian{Y_p}{\parallelTransport{q}{p} Z_q}
\end{equation}
for all $Y_p \in \tangentSpace{p}$ and all $Z_q \in \tangentSpace{q}$.
Moreover, we have
\begin{equation}
	\label{eq:logarithm:swapping-the-points}
	\parallelTransport{p}{q} \logarithm{p} q
	=
	- \logarithm{q} p
	.
\end{equation}

\begin{definition}
	A function $f \colon \cM \to \eR \coloneqq \bbR \cup \{+\infty\}$ is said to be \emph{lower semi-continuous} at $p \in \cM$ if
	\begin{equation*}
		\liminf_{q \to p} f(q)
		\ge
		f(p)
		.
	\end{equation*}
\end{definition}

\begin{definition}
	\label{definition:geodesic-convexity}
	Let $\cU \subseteq \cM$ be a strongly convex set.
	A function $f \colon \cM \to \eR$ is said to be \emph{geodesically convex on~$\cU$} if for any $p, q \in \cU$,
	\begin{equation}
		\label{eq:geodesic-convexity}
		(f \circ \geodesic{p}{q})(t)
		\le
		(1 -t) f(p) + t f(q)
		\quad
		\text{holds for all }
		t \in [0,1]
		.
	\end{equation}
\end{definition}

A function as in \zcref{definition:geodesic-convexity} can be obtained, for instance, by starting with a strongly convex set $\cU$ and a function $f \colon \cU \to \bbR$ satisfying \zcref{eq:geodesic-convexity} for any $p, q \in \cU$, and by extending it to $+\infty$ outside of~$\cU$.
Then $f$ will be geodesically convex on $\dom f \coloneqq \setDef[auto]{p \in \cM}{f(p) \ne +\infty} = \cU$.

Let $f \colon \cM \to \eR$ be a function such that $\dom f$ is strongly convex and $f$ is a geodesically convex function on $\dom f$.
The \emph{subdifferential} of $f$ at $p \in \dom f$ is defined as
\begin{equation*}
	\partial f(p)
	\coloneqq
	\setDef[big]{X_p \in \tangentSpace{p}}{f(q) \ge f(p)
		+
	\riemannian{X_p}{\logarithm{p} q} \text{ for all } q \in \dom f}
	.
\end{equation*}
For any $\varepsilon > 0$, the $\varepsilon$-subdifferential of $f$ at $p \in \dom f$ is defined as
\begin{equation*}
	\partial_\varepsilon f(p)
	=
	\setDef[big]{X_p \in \tangentSpace{p}}{f(q) \ge f(p)
		+
	\riemannian{X_p}{\logarithm{p} q} - \varepsilon \text{ for all } q \in \dom f}
	.
\end{equation*}
Clearly, the inclusion
\begin{equation*}
	\partial f(p)
	\subseteq
	\partial_\varepsilon f(p)
\end{equation*}
holds for all $\varepsilon > 0$.

\section{Approximation of the \texorpdfstring{$\varepsilon$}{epsilon}-Subdifferential}
\label{section:approximation-of-the-epsilon-subdifferential}

We begin by specifying the standing assumptions for the remainder of the paper.
\begin{assumption}[standing assumptions] \hfill
	\label{assumption:standing-assumptions}
	\begin{enumerate}
		\item
			$\cM$ is a complete Riemannian manifold of finite dimension.

		\item
			$f \colon \cM \to \eR$ is a lower semi-continuous function.

		\item
			$\dom f \subseteq \cM$ is a strongly convex set with nonempty interior in~$\cM$, and $f$ is a geodesically convex function on $\dom f$.

		\item
      There exist two numbers $\omega, \Omega \in \bbR$ such that the sectional curvature satisfies $\omega \le K_p(X_p, Y_p) \le \Omega$ for all $p \in \dom f$ and all linearly independent unit vectors $X_p, Y_p \in \tangentSpace{p}$.

		\item
			In case $\omega < 0$ or $\Omega > 0$, $\dom f$ is bounded and $\delta \coloneqq \diam(\dom f)$ denotes its diameter.
			In case $\Omega > 0$, we assume $\delta < \frac{\pi}{\sqrt{\Omega}}$.
			If the manifold is flat ($\omega = \Omega = 0$), then $\dom f$ is allowed to be unbounded.
	\end{enumerate}
\end{assumption}
Notice that the third assumption in particular implies that the logarithmic map is well-defined on $\dom f$.
The last assumption is due to the fact that in case of $\Omega > 0$, sets of diameter larger than $\frac{\pi}{\sqrt{\Omega}}$ cannot be strongly convex; see, \eg, \cite[Proposition~4.2~(iii)]{LiYao:2012:1}, \cite[Proposition~4.1]{WangLiYao:2016:1} or \cite[Proposition~II.1.4]{BridsonHaefliger:1999:1}.

As usual in bundle methods see, \eg, \cite[Chapter~10]{BonnansGilbertLemarechalSagastizabal:2006:1}, we aim to approximate the objective function through a sequence of models $\sequence{\varphi}{k} \colon \dom f \to \bbR$, $k \in \bbN_0$, defined by
\begin{equation*}
	\sequence{\varphi}{k}(q)
	\coloneqq
	\max \setDef[big]{f(\sequence{p}{j}) + \riemannian{X_{\sequence{p}{j}}}{\logarithm{\sequence{p}{j}} q}}{j = 0, \ldots, k}
	,
	\quad
	q \in \dom f
	,
\end{equation*}
where $X_{\sequence{p}{j}} \in \partial f(\sequence{p}{j})$ is a \emph{subgradient} of $f$ at $\sequence{p}{j} \in \cM$ for all $j = 0, \ldots, k \in \bbN_0$.
Define now the \emph{linearization errors} as
\begin{equation}
	\label{eq:new-lin-errors}
	\sequence{e_j}{k}
	\coloneqq
	f(\sequence{p}{k}) - f(\sequence{p}{j}) - \riemannian{X_{\sequence{p}{j}}}{\logarithm{\sequence{p}{j}} \sequence{p}{k}}
\end{equation}
and note that $\sequence{e_j}{k} \ge 0$ holds for all $j = 0, \ldots, k$ by the definition of subgradient.
We can therefore rewrite
\begin{equation} \label{eq:model}
	\sequence{\varphi}{k}(q)
	=
	f(\sequence{p}{k})
	+
	\max \setDef[big]{- \sequence{e_j}{k} + \riemannian{X_{\sequence{p}{j}}}{\logarithm{\sequence{p}{j}} q - \logarithm{\sequence{p}{j}} \sequence{p}{k}}}{j = 0, \ldots, k}
\end{equation}
for $	q \in \dom f$.

The linearization errors defined in \zcref{eq:new-lin-errors} and the representation of the model in \zcref{eq:model} are the exact analogs of their counterparts for convex bundle methods in vector spaces.
However, it turns out that, contrary to the vector space situation, the linearization errors cannot be used directly to define approximations of the $\varepsilon$-subdifferential of~$f$.
This is because it is not clear a priori that such a choice of linearization errors ensures that a subgradient at an initial point~$p$ can be parallely transported to a different point~$q$ and remain an element of the $\varepsilon$-subdifferential at~$q$, as is the case when dealing with vector spaces.
This is an effect caused by curvature.

For $p \in \dom f$ we introduce the following quantities,
\begin{equation}
	\label{eq:linearization-curvature-remainder}
	\sequence{r_j}{k}(p)
	\coloneqq
	\varrho_{\omega, \Omega}(\dist(\tilde p, \, p))
	\,
	\riemanniannorm{X_{\sequence{p}{j}}}
	\,
	\riemanniannorm{\logarithm{\sequence{p}{j}}\sequence{p}{k}}
	,
\end{equation}
where $\tilde p$ lies on the minimal geodesic arc that connects $\sequence{p}{j}$ to $\sequence{p}{k}$,
\begin{equation}
	\label{eq:rho-definition}
	\varrho_{\omega, \Omega}(s)
	\coloneqq
	\max \paren[big]\{\}{%
		\zeta_{1, \omega}(s) - 1
		,
		\;
		1 - \zeta_{2, \Omega}(s)
	}
	,
\end{equation}
and
\begin{equation}
	\label{eq:hessian-eigs}
	\begin{aligned}
		\zeta_{1, \omega}(s)
		&
		\coloneqq
		\begin{cases}
			1
			&
			\text{if }
			\omega \ge 0
			,
			\\
			\sqrt{- \omega} \, s \coth(\sqrt{- \omega} \, s)
			\quad
			&
			\text{if }
			\omega < 0
			,
		\end{cases}
		\\
		\zeta_{2, \Omega}(s)
		&
		\coloneqq
		\begin{cases}
			1
			&
			\text{if }
			\Omega \le 0
			,
			\\
			\sqrt{\Omega} \, s \cot(\sqrt{\Omega} \, s)
			\quad
			&
			\text{if }
			\Omega > 0
			,
		\end{cases}
	\end{aligned}
\end{equation}
for $s \in \bbR$.
We refer the reader to \zcref{figure:hessian-eigs} for an illustration.
Observe that
\begin{equation}
	\label{eq:remainder-upper-bound}
	\sequence{r_j}{k}(p)
	\le
	\varrho_{\omega, \Omega}(\delta)
	\,
	\riemanniannorm{X_{\sequence{p}{j}}}
	\,
	\riemanniannorm{\logarithm{\sequence{p}{j}}\sequence{p}{k}}
\end{equation}
holds for all $p \in \dom f$.
For more details, see \zcref{section:curvature-bound}.

As the curvature bounds~$\omega$ and~$\Omega$, as well as the diameter~$\delta$ are fixed throughout the manuscript, we will use $\varrho$ as a shorthand notation for $\varrho_{\omega, \Omega}(\delta)$.
Then $\sequence{r_j}{k}(p)$ simplifies to
\begin{equation}
	\label{eq:fixed-remainder}
	\sequence{r_j}{k}
	\coloneqq
	\varrho
	\,
	\riemanniannorm{X_{\sequence{p}{j}}}
	\,
	\riemanniannorm{\logarithm{\sequence{p}{j}}\sequence{p}{k}}
	.
\end{equation}
The curvature-dependent quantities $\sequence{e_j}{k} + \sequence{r_j}{k}$ allow for an approximation of the $\varepsilon$-subdifferential of~$f$ at a given point on the manifold using a convex combination of parallely transported subgradients from other points.
This approximation depends on bounds that involve the sectional curvature of the manifold.
\begin{figure}
	\centering
	\begin{tikzpicture}[
		> = stealth,
		declare function = {
			zeta1(\x,\omega) = sqrt(-\omega) * x / tanh(sqrt(-\omega) * x) * (abs(x) > 0.1) + 1 * (abs(x) <= 0.1);
		},
		baseline
		]
		\begin{axis}[
			axis x line = center,
			axis y line = center,
			width = 0.49\textwidth,
			xlabel = $\delta$,
			xlabel style = {at = {(axis cs:2.25, 0)}, anchor = north west},
			enlarge x limits = 0.15,
			xmin = 0,
			xmax = 2,
			ymin = -0.33,
			ymax = 3.3,
			ytick = {1, 2, 3},
			xtick = {1, 2},
			cycle list name = TolLight,
			legend style = {at = {(0.98, 0.12)}, anchor = south east, font = \footnotesize},
			]
			\edef\lstomega{-1.0, -2.0, -3.0}
			\pgfplotsinvokeforeach{0, 1, 2}{
				\pgfmathsetmacro{\myomega}{{\lstomega}[#1]}
				\edef\temp{%
					\noexpand\addplot+ [domain = 0:2, samples = 101, mark = none, thick]{zeta1(x,\myomega)};
					\noexpand\addlegendentry{$\omega = \myomega$};
				}
				\temp
			}
		\end{axis}
	\end{tikzpicture}
	\hspace{.005\textwidth}
	\begin{tikzpicture}[
		>=stealth,
		declare function = {
			zeta2(\x,\Omega) = 1 - \Omega/3 * x^2 - \Omega^2/45 * x^4 - 2*\Omega^3/945 * x^6;
		},
		baseline
		]
		\begin{axis}[
			axis x line = center,
			axis y line = center,
			width = 0.49\textwidth,
			xlabel = $\delta$,
			xlabel style = {at = {(axis cs:1.9, 0)}, anchor = north west},
			enlarge x limits = 0.15,
			xmin = 0,
			xmax = 1.7,
			ymin = -0.1,
			ymax = 1.1,
			xtick = {pi/(2*sqrt(3)), pi/(2*sqrt(2)), pi/(2)},
			xticklabels = {$\!\!\!\frac{\pi}{2 \sqrt{3}}$, $\,\,\,\frac{\pi}{2 \sqrt{2}}$, $\frac{\pi}{2}$},
			ytick = {0.5, 1},
			yticklabels = {$\frac{1}{2}$, $1$},
			xticklabel style = {yshift = 1.05mm},
			cycle list name = TolLight,
			legend style = {at = {(0.98, 0.98)}, anchor = north east, font = \footnotesize},
			]
			\edef\lstOmega{1.0, 2.0, 3.0}
			\pgfplotsinvokeforeach{0, 1, 2}{
				\pgfmathsetmacro{\myOmega}{{\lstOmega}[#1]}
				\pgfmathsetmacro{\myMax}{pi / (2*sqrt(\myOmega))}
				\edef\temp{%
					\noexpand\addplot+ [domain = 0:\myMax, samples = 101, mark = none, thick]{zeta2(x,\myOmega)};
					\noexpand\addlegendentry{$\Omega = \myOmega$};
				}
				\temp
			}
		\end{axis}
	\end{tikzpicture}
	\caption{Illustration of the functions $\zeta_{1, \omega}$ (left) and $\zeta_{2, \Omega}$ (right) defined in \eqref{eq:hessian-eigs}.}
	\label{figure:hessian-eigs}
\end{figure}

We are now in a position to prove a result parallel to \cite[Theorem~6.68]{GeigerKanzow:2002:1}, which shows that we obtain an inner approximation of the $\varepsilon$-subdifferential, using convex combinations of parallely transported subgradients.
In what follows, we will occasionally use the symbol $\Delta_k$ to denote the probability simplex in $\bbR^{k+1}$, \ie, $\Delta_k \coloneqq \setDef[big]{\lambda \in \bbR^{k+1}}{\sum_{j=0}^k \lambda_j = 1, \; \lambda_j \ge 0, \; j = 0, \ldots, k}$.
\begin{theorem}
	\label{theorem:eps-approx}
	Suppose that $k \in \bbN_0$, $\sequence{p}{j} \in \dom f$ and $X_{\sequence{p}{j}} \in \partial f(\sequence{p}{j})$ holds for all $j = 0, \ldots, k$.
	Then, for any $\varepsilon > 0$, the subset of $\tangentSpace{\sequence{p}{k}}$ defined by
	\begin{equation}
		\label{eq:eps-approx}
		\sequence{G_\varepsilon}{k}
		\coloneqq
		\setDef[Big]{
			\sum_{j=0}^k
			\lambda_j
			\,
			\parallelTransport{\sequence{p}{j}}{\sequence{p}{k}}
			X_{\sequence{p}{j}}
		}{
			\sum_{j=0}^k
			\lambda_j
			\,
			\paren[big](){%
				\sequence{e_j}{k}
				+
				\sequence{r_j}{k}
			}
			\le
			\varepsilon
			,
			\;
			\lambda \in \Delta_k
		}
	\end{equation}
	satisfies
	\begin{equation*}
		\sequence{G_\varepsilon}{k}
		\subseteq
		\partial_\varepsilon f(\sequence{p}{k})
		.
	\end{equation*}
\end{theorem}
\begin{proof}
	Let $p \in \dom f$, then
	\begin{align*}
		\MoveEqLeft
		\riemannian[big]{X_{\sequence{p}{j}}}{\logarithm{\sequence{p}{j}} p - \logarithm{\sequence{p}{j}} \sequence{p}{k}}
		\\
		&
		=
		\riemannian[big]{X_{\sequence{p}{j}}}{\logarithm{\sequence{p}{j}} p}
		-
		\riemannian[big]{X_{\sequence{p}{j}}}{\logarithm{\sequence{p}{j}} \sequence{p}{k}}
		\\
		&
		\le
		f(p)
		-
		f(\sequence{p}{j})
		-
		\riemannian[big]{X_{\sequence{p}{j}}}{\logarithm{\sequence{p}{j}} \sequence{p}{k}}
		\\
		&
		=
		f(p)
		-
		f(\sequence{p}{k})
		+
		\sequence{e_j}{k}
		,
	\end{align*}
	where we used the subdifferential inequality at $\sequence{p}{j}$ and the definition of linearization errors given in \zcref{eq:new-lin-errors}.
	Therefore
	\begin{equation*}
		f(p)
		\ge
		f(\sequence{p}{k})
		+
		\riemannian[big]{X_{\sequence{p}{j}}}{
			\logarithm{\sequence{p}{j}} p
			-
			\logarithm{\sequence{p}{j}} \sequence{p}{k}
		}
		-
		\sequence{e_j}{k}
		\quad
		\text{for all }
		p \in \dom f
		.
	\end{equation*}
	Hence, we see that the model $\sequence{\varphi}{k}$ approximates the objective function $f$ from below on~$\dom f$ due to the linearization errors \zcref{eq:new-lin-errors} satisfying $\sequence{e_j}{k} \ge 0$.

	By adding and subtracting $\riemannian[big]{\parallelTransport{\sequence{p}{j}}{\sequence{p}{k}} X_{\sequence{p}{j}}}{\logarithm{\sequence{p}{k}} p}$ to the right-hand side of the previous inequality, we get
	\begin{align*}
		f(p)
		&
		\ge
		f(\sequence{p}{k})
		+
		\riemannian[big]{
			\parallelTransport{\sequence{p}{j}}{\sequence{p}{k}}
			X_{\sequence{p}{j}}
		}{
			\logarithm{\sequence{p}{k}} p
		}
		\\
		&
		\quad
		+
		\riemannian[big]{X_{\sequence{p}{j}}}{
			\logarithm{\sequence{p}{j}} p
			-
			\parallelTransport{\sequence{p}{k}}{\sequence{p}{j}}
			\logarithm{\sequence{p}{k}} p
			-
			\logarithm{\sequence{p}{j}} \sequence{p}{k}
		}
		-
		\sequence{e_j}{k}
		,
	\end{align*}
	for all $p \in \dom f$.
	Following an argument presented in \cite[Appendix~C]{AlimisisOrvietoBecigneulLucchi:2021:1}, see also~\zcref{section:curvature-bound}, we can write
	\begin{equation}
		\label{ineq:almost-subgradient}
		f(p)
		\ge
		f(\sequence{p}{k})
		+
		\riemannian[big]{
			\parallelTransport{\sequence{p}{j}}{\sequence{p}{k}}
			X_{\sequence{p}{j}}
		}{
			\logarithm{\sequence{p}{k}} p
		}
		-
		\sequence{r_j}{k}
		-
		\sequence{e_j}{k}
		,
	\end{equation}
	for all $p \in \dom f$, with $\sequence{r_j}{k}$ from \zcref{eq:fixed-remainder}.

	Let $\lambda \in \Delta_k$, multiply \zcref{ineq:almost-subgradient} by $\lambda_j$ and sum up to obtain
	\begin{equation}
		\label{eq:epsilon-subdifferential-inequality}
		f(p)
		\ge
		f(\sequence{p}{k})
		+
		\sum_{j=0}^k \lambda_j \, \riemannian[big]{\parallelTransport{\sequence{p}{j}}{\sequence{p}{k}} X_{\sequence{p}{j}}}{\logarithm{\sequence{p}{k}} p}
		-
		\sum_{j=0}^k \lambda_j
		\,
		\paren[big](){%
			\sequence{e_j}{k}
			+
			\sequence{r_j}{k}
		}
		,
	\end{equation}
	for all $p \in \dom f$.
	Hence, by bi-linearity of the Riemannian metric, we can take the sum and multiplication by $\lambda_j$ inside the inner product.
	Therefore
	\begin{equation*}
		\sum_{j=0}^k \lambda_j \, \parallelTransport{\sequence{p}{j}}{\sequence{p}{k}} X_{\sequence{p}{j}}
		\in
		\partial_\varepsilon f(\sequence{p}{k})
	\end{equation*}
	holds if and only if
	\begin{equation*}
		\sum_{j=0}^k \lambda_j \,
		\paren[big](){%
			\sequence{e_j}{k}
			+
			\sequence{r_j}{k}
		}
		\le
		\varepsilon
		,
	\end{equation*}
	which concludes the argument.
\end{proof}

Note that the set $\sequence{G_\varepsilon}{k}$ is a subset of the inner-product space $\tangentSpace{\sequence{p}{k}}$, and it is nonempty, convex, and compact.
The orthogonal projection $\proj \colon \tangentSpace{\sequence{p}{k}} \to \sequence{G_\varepsilon}{k}$ of $Y_{\sequence{p}{k}}$, given by
\begin{equation*}
	\riemanniannorm{Y_{\sequence{p}{k}} - \proj{\sequence{G_\varepsilon}{k}}(Y_{\sequence{p}{k}})}
	\le
	\riemanniannorm{Y_{\sequence{p}{k}} - Z_{\sequence{p}{k}}}
	\quad
	\text{for all }
	Z_{\sequence{p}{k}}
	\in
	\sequence{G_\varepsilon}{k}
	,
\end{equation*}
is thus well-defined; see for example~\cite[Lemma~2.17]{GeigerKanzow:2002:1} for a proof of this fact.
We can therefore define
\begin{equation}
	\label{eq:search-dir-proj}
	\sequence{g}{k}
	\coloneqq
	\proj{\sequence{G_\varepsilon}{k}}(0_{\sequence{p}{k}})
	.
\end{equation}
The negative of this tangent vector will serve as search direction in our bundle method.

\zcref{eq:search-dir-proj} means $\riemanniannorm{\sequence{g}{k}} \le \riemanniannorm{Z_{\sequence{p}{k}}}$ for all $Z_{\sequence{p}{k}} \in \sequence{G_\varepsilon}{k}$, or equivalently
\begin{equation*}
	\sequence{g}{k}
	\coloneqq
	\argmin_{Z_{\sequence{p}{k}} \in \sequence{G_\varepsilon}{k}} \riemanniannorm{Z_{\sequence{p}{k}}}
	.
\end{equation*}
Using the characterization of $\sequence{G_\varepsilon}{k}$ from \zcref{theorem:eps-approx}, we can therefore write
\begin{equation}
	\label{eq:search-dir-sum}
	\sequence{g}{k}
	=
	\sum_{j=0}^k \sequence{\lambda_j}{k} \, \parallelTransport{\sequence{p}{j}}{\sequence{p}{k}} X_{\sequence{p}{j}}
	,
\end{equation}
where the tuple $\sequence{\lambda}{k} = (\sequence{\lambda}{k}_{0}, \ldots, \sequence{\lambda_k}{k})$ is a solution to the following convex quadratic program,
\begin{equation}
	\label{problem:P1}
  \zcsetup{reftype=problem}
	\tag{P1}
	\begin{aligned}
		\argmin_{\lambda \in \bbR^{k+1}}
		\quad
		&
		\frac{1}{2} \riemanniannorm[Big]{\sum_{j=0}^k \lambda_j \, \parallelTransport{\sequence{p}{j}}{\sequence{p}{k}} X_{\sequence{p}{j}}}^2
		\\
		\text{\suto}
		\quad
		&
		\sum_{j=0}^k \lambda_j \,
		\paren[big](){%
			\sequence{e_j}{k}
			+
			\sequence{r_j}{k}
		}
		\le
		\varepsilon
		,
		\quad
		\text{and}
    \quad
		\lambda \in \Delta_k
		.
	\end{aligned}
\end{equation}

Closely related to \zcref{problem:P1} is another quadratic program, in which the weighted sum of the components of $\sequence{\lambda}{k}$ is not constrained but rather penalized:
\begin{equation}
  \zcsetup{reftype=problem}
	\label{problem:P2}
	\tag{P2}
	\begin{aligned}
		\argmin_{\lambda \in \bbR^{k+1}}
		\quad
		&
		\frac{1}{2} \riemanniannorm[Big]{
			\sum_{j=0}^k
			\lambda_j
			\,
			\parallelTransport{\sequence{p}{j}}{\sequence{p}{k}} X_{\sequence{p}{j}}
		}^2
		+
		\sum_{j=0}^k
		\lambda_j
		\,
		\paren[big](){%
			\sequence{e_j}{k}
			+
			\sequence{r_j}{k}
		}
		\\
		\text{\suto}
		\quad
		&
		\lambda \in \Delta_k
		.
	\end{aligned}
\end{equation}
The relation between \zcref{problem:P1} and \zcref{problem:P2} can be studied through the following result, the proof of which is analogous to \cite[Lemma~6.70]{GeigerKanzow:2002:1}.

\begin{lemma}
	\label{lemma:kktp2}
	The vector $\sequence{\lambda}{k} = (\sequence{\lambda_0}{k}, \ldots, \sequence{\lambda_k}{k}) \in \bbR^{k+1}$ is a solution to \zcref{problem:P2} if and only if there exist $\sequence{g}{k} \in \tangentSpace{\sequence{p}{k}}$ and a Lagrange multiplier $\sequence{\xi}{k} \in \bbR$ such that the triple $(\sequence{\lambda}{k}, \sequence{g}{k}, \sequence{\xi}{k})$ fulfills the following KKT~conditions
	\begin{align*}
		g
		=
		\sum_{j=0}^k \lambda_j \, \parallelTransport{\sequence{p}{j}}{\sequence{p}{k}} X_{\sequence{p}{j}}
		,
    \quad
    &
		\\
		\sum_{j=0}^k \lambda_j
		=
		1
    \quad
		\text{and}
    \quad
		\lambda_j
		\ge
		0
    \quad
    &
		\text{for all }
		j = 0, \ldots, k
		,
		\\
		\riemannian{g}{\parallelTransport{\sequence{p}{j}}{\sequence{p}{k}} X_{\sequence{p}{j}}}
		+
		\sequence{e_j}{k}
		+
		\sequence{r_j}{k}
		+
		\xi
		\ge
		0
		\quad
    &
		\text{for all }
		j = 0, \ldots, k
		,
		\\
		\lambda_j
		\paren[big](){%
			\riemannian{g}{
				\parallelTransport{\sequence{p}{j}}{\sequence{p}{k}}
			X_{\sequence{p}{j}}}
			+
			\sequence{e_j}{k}
			+
			\sequence{r_j}{k}
			+
			\xi
		}
		=
		0
		\quad
    &
		\text{for all }
		j = 0, \ldots, k
		.
	\end{align*}
\end{lemma}

Let now $\sequence{\lambda}{k} = (\sequence{\lambda_0}{k}, \ldots, \sequence{\lambda_k}{k}) \in \bbR^{k+1}$ be a solution to \zcref{problem:P2} and choose $\varepsilon \coloneq \sum_{j = 0}^{k} \sequence{\lambda_j}{k} \,\paren[big](){\sequence{e_j}{k} + \sequence{r_j}{k}}$ analogously to \cite[p.389]{GeigerKanzow:2002:1}.
We find that the KKT conditions of \zcref{problem:P1}, given by
	\begin{align*}
		g
		=
		\sum_{j=0}^k \lambda_j \, \parallelTransport{\sequence{p}{j}}{\sequence{p}{k}}
    X_{\sequence{p}{j}}
    &
		,
		\\
		\sum_{j=0}^k \lambda_j
		=
		1
    \quad
		\text{and}
    \quad
		\lambda_j
		\ge
		0
    &
    \quad
		\text{for all }
		j = 0, \ldots, k
		,
		\\
		\riemannian{g}{\parallelTransport{\sequence{p}{j}}{\sequence{p}{k}} X_{\sequence{p}{j}}}
		+
    \zeta (
      \sequence{e_j}{k}
      +
      \sequence{r_j}{k}
    )
		+
		\xi
		\ge
		0
    &
		\quad
		\text{for all }
		j = 0, \ldots, k
		,
		\\
		\lambda_j
		\paren[big](){%
			\riemannian{g}{
				\parallelTransport{\sequence{p}{j}}{\sequence{p}{k}}
			X_{\sequence{p}{j}}}
			+
      \zeta (
        \sequence{e_j}{k}
        +
        \sequence{r_j}{k}
      )
			+
			\xi
		}
		=
		0
    &
		\quad
		\text{for all }
		j = 0, \ldots, k
    ,
    \\
    \zeta \big(
     \sum_{j = 0}^{k}
     \sequence{\lambda_j}{k} \,\paren[big](){
      \sequence{e_j}{k}
      +
      \sequence{r_j}{k}
    }
    -
    \varepsilon
    \big)
    =
    0
    \quad
    \text{and}
    \quad
    \zeta
    \ge
    0
    &
    ,
    \\
    \sum_{j = 0}^{k}
    \sequence{\lambda_j}{k} \,\paren[big](){
      \sequence{e_j}{k}
      +
      \sequence{r_j}{k}
    }
    \le
    \varepsilon
    &
		,
	\end{align*}
are also satisfied by solutions of \zcref{problem:P2}.
Thus, solutions of \zcref{problem:P2} are solutions to \zcref{problem:P1} for this choice of $\varepsilon$.
Lastly, we can compute the unique Lagrange multiplier of \zcref{problem:P2} from $\sequence{\lambda}{k}$ as follows:
\begin{align*}
	\sequence{\xi}{k}
	&
	=
	\sequence{\xi}{k} \sum_{j=0}^k \sequence{\lambda_j}{k}
	=
	\sum_{j=0}^k
	\sequence{\lambda_j}{k}
	\,
	\paren[big](){%
		-
		\riemannian[big]{\sequence{g}{k}}{
			\parallelTransport{\sequence{p}{j}}{\sequence{p}{k}}
			X_{\sequence{p}{j}}
		}
		-
		\sequence{e_j}{k}
		-
		\sequence{r_j}{k}
	}
	\\
	&
	=
	-
	\riemannian[big]{\sequence{g}{k}}{
		\sum_{j=0}^k
		\sequence{\lambda_j}{k}
		\,
		\parallelTransport{\sequence{p}{j}}{\sequence{p}{k}}
		X_{\sequence{p}{j}}
	}
	-
	\sum_{j=0}^k
	\sequence{\lambda_j}{k}
	\,
	\sequence{e_j}{k}
	-
	\sum_{j=0}^k
	\sequence{\lambda_j}{k}
	\,
	\sequence{r_j}{k}
	\\
	&
	=
	-
	\riemanniannorm{\sequence{g}{k}}^2
	-
	\sum_{j=0}^k \sequence{\lambda_j}{k} \, \sequence{e_j}{k}
	-
	\sum_{j=0}^k
	\sequence{\lambda_j}{k}
	\,
	\sequence{r_j}{k}
	.
\end{align*}

\begin{remark}
	\label{remark:retraction-vector-transport-generalization}
	It is possible to utilize more general tools than the exponential map and parallel transport in \zcref{theorem:eps-approx}.
	We introduce two ways this could be achieved.

	First, let $\retract{p}{X_p}$ be a retraction (see, \eg, \cite[Definition~4.1.1]{AbsilMahonySepulchre:2008:1}) of $X_p \in \tangentSpace{p}$ on $\cM$.
  Let $\vectorTransportDir{p}{Y_p} X_p$ be a continuous vector transport (see, \eg, \cite[Definition~8.1.1]{AbsilMahonySepulchre:2008:1}) of $X_p \in \tangentSpace{p}$ from $p \in \cM$ to the point $\retract{p}Y_p \in \cM$, and let $\vectorTransportDir{p}{Y_p}[\retractionSymbol] X = \d_{Y_p} \retract{p}[X_p]$ be the differential of the retraction.

	As in \cite[Assumption~1.3]{KasaiSatoMishra:2018:1}, or, alternatively, \cite[(2.5)]{HuangGallivanAbsil:2015:1}, assume that there exists a constant $C_0 > 0$ such that $\vectorTransportSymbol$ satisfies
	\begin{equation*}
		\riemanniannorm{
			\vectorTransportDir{p}{Y_p}
			-
			\vectorTransportDir{p}{Y_p}[\retractionSymbol]
		}
		\le
		C_0
		\riemanniannorm{Y_p}
		,
	\end{equation*}
	for all $p, q \in \dom f$ with $\retract{p}{Y_p} \eqcolon q$.
  We remark that this inequality is satisfied for all such $p, q$ as long as the vector transport $\vectorTransport{p}{q} \coloneqq \vectorTransportDir{p}{Y_p}$ is continuous.
	Then, as in \cite[Lemma~3.7]{KasaiSatoMishra:2018:1}, or, alternatively, \cite[Lemma~3.5]{HuangGallivanAbsil:2015:1}, there exists a constant $C_t > 0$ such that, for all $p, q \in \dom f$
	\begin{equation*}
		\riemanniannorm{
			\vectorTransport{p}{q} X_p
			-
			\parallelTransport{p}{q} X_p
		}
		\le
		C_t
		\riemanniannorm{Y_p}
		\riemanniannorm{X_p}
		,
	\end{equation*}
	where $X_p, Y_p \in \tangentSpace{p}$ and $\retract{p}{Y_p} = q$.
	Similarly, assume as in \cite[Assumption~1.6]{KasaiSatoMishra:2018:1} that there exists a constant $C_r > 0$ such that
	\begin{equation}
		\label{eq:inverse-retraction-approximation}
		\riemanniannorm{
			\logarithm{p} q
			-
			\inverseRetract{p} q
		}
		\le
		C_r
		\,
		\riemanniannorm{\inverseRetract{p} q}^2
		,
	\end{equation}
	for all $p, q \in \dom f$.
	It is also necessary to assume that the retraction is invertible in~$\dom f$.
	Note that \zcref{eq:inverse-retraction-approximation} is quite a natural assumption as retractions are first order approximations of the exponential map.
	Let now
	\begin{equation}
		\label{eq:inverse-retraction-linearization-errors}
		\sequence{\tilde e_j}{k}
		\coloneqq
		f(\sequence{p}{k})
		-
		f(\sequence{q}{j})
		-
		\riemannian{
			X_{\sequence{q}{j}}
		}{
			\inverseRetract{\sequence{q}{j}} \sequence{p}{k}
		}
		,
	\end{equation}
	and
	\begin{equation}
		\label{eq:inverse-retraction-curvature-remainder}
		\sequence{\tilde r_j}{k}
		\coloneqq
		\varrho \, \riemanniannorm[big]{X_{\sequence{q}{j}}} \, \riemanniannorm[big]{\inverseRetract{\sequence{q}{j}} \sequence{p}{k+1}}
		.
	\end{equation}
	Combining these facts with \zcref{eq:epsilon-subdifferential-inequality} we get
	\begin{equation*}
		\sum_{j=0}^k \lambda_j
		\,
		\vectorTransport{\sequence{q}{j}}{\sequence{p}{k}}
		X_{\sequence{q}{j}}
		\in
		\partial_{\varepsilon_e + \varepsilon_r + \varepsilon_t} f(\sequence{p}{k})
	\end{equation*}
	if and only if
	$
	\displaystyle\sum_{j=0}^k \lambda_j
	\paren[big](){%
		\sequence{\tilde e_j}{k}
		+
		\sequence{\tilde r_j}{k}
	}
	\le
	\varepsilon_e
	,
	$
	as before, and additionally
	\begin{equation*}
		\begin{aligned}
			&
			\sum_{j=0}^k \lambda_j
			\,
			C_r \, (1 + \varrho)
			\,
			\riemanniannorm{X_{\sequence{q}{j}}}
			\,
			\riemanniannorm{\inverseRetract{\sequence{q}{j}} \sequence{p}{k}}^2
			\le
			\varepsilon_r
			\\
			\text{and}
			\quad
			&
			\sum_{j=0}^k \lambda_j
			\,
			C_t
			\,
			\delta
			\,
			\riemanniannorm{X_{\sequence{q}{j}}}
			\riemanniannorm{
				\logarithm{\sequence{q}{j}} \sequence{p}{k}
			}
			\le
			\varepsilon_t
			,
		\end{aligned}
	\end{equation*}
	for all $k \in \bbN$ and arbitrary $\varepsilon_e, \varepsilon_r, \varepsilon_t > 0$.
	Hence, for $\varepsilon \coloneqq \varepsilon_e + \varepsilon_r + \varepsilon_t$, the quality of the inner approximation of the $\varepsilon$-subdifferential of $f$ at $\sequence{p}{k}$ depends on $C_r$ and $C_t$ as well.

	As a second possible approach, one may rephrase all definitions related to geodesic convexity using \emph{retraction convexity} instead, as introduced, \eg, in \cite[Definition~2]{HuangWei:2021:1}.
	This would entail replacing geodesic arcs $\geodesic{p}{q}(t) = \exponential{p}(t \logarithm{p}q)$ with retraction curves $c_{p, q}(t) = \retract{p}(t \, \inverseRetract{p}q)$.
	Most notably, the definition of $\varepsilon$-subdifferential would then include the inverse of the chosen retraction, which is thus evidently required to be invertible in $\dom f$.
	Furthermore, for the derivation in \zcref{section:curvature-bound} to hold in this case, it is necessary to use a vector transport that is parallel along the aforementioned retraction curves.
\end{remark}

\section{Presentation of the Method}
\label{section:algorithm}

In this section we present the proposed Riemannian convex bundle method as \zcref{algorithm:RCBM}.
The notation is close to the one in \cite[Chapter~6.7]{GeigerKanzow:2002:1}.

After an initialization step, in \zcref{line:subprob} the method computes a solution $\sequence{\lambda}{k}$ of the convex, quadratic \zcref{problem:P2} formed using the subgradients currently held in the bundle, indexed by $\sequence{J}{k}$.
This vector of convex coefficients is then used in \zcref{line:search-direction} to procure a search direction $\sequence{d}{k}$, as well as the numbers $\sequence{\varepsilon}{k}$ and $\sequence{\sigma}{k}$ that quantify how well the subdifferential is approximated at the $k$-th step.
The tolerance given by the user prescribes whether such approximation is satisfactory, and whether the zero vector is close to being an element thereof, in which case the algorithm is stopped.
A new candidate point is constructed by moving from the last serious iterate in the new search direction with a suitable step size (\zcref{line:backtracking-loop} to \zcref{line:end-backtracking-loop}).
To assess whether the new candidate is a serious iterate, a check is carried out (\zcref{line:descent-test}) to verify whether the objective decreased enough.
If this check is successful, the candidate is labeled as the new serious iterate.
Otherwise, the algorithm performs a null step and looks for a new candidate via computing a smaller step size until a certain geometric condition is met.
The new candidate and a subgradient of the objective sampled at this point are appended to the bundle irrespectively of whether a serious or a null step is being performed.
Finally, the algorithm purges unused subgradients from the bundle and updates the linearization errors and curvature-dependent remainder terms (\zcref{line:positive-coefficients} to \zcref{line:curvature-remainder}).

\begin{algorithm}[htp]
	\caption{Riemannian Convex Bundle Method (RCBM)}
	\label{algorithm:RCBM}
	\begin{algorithmic}[1]
		\Require
		$f$,
		oracle for~$\partial f$,
		initial point $\sequence{p}{0} \in \interior \dom f$
		\Require
		parameters
		$\varrho$,
		$m \in~(0, 1)$,
		$\beta \in (0,1)$,
		$\tol > 0$
		\Ensure
		last serious iterate $\sequence{p}{k_\ast}$ for some $k_\ast \in \bbN_0$
		\State
		Set $\sequence{q}{0} \coloneqq \sequence{p}{0}$, $X_{\sequence{q}{0}} \in \partial f(\sequence{p}{0})$, $\sequence{g}{0} \coloneqq X_{\sequence{q}{0}}$, $\sequence{\varepsilon}{0} \coloneqq 0$, $\sequence{e_0}{0} \coloneqq 0$, $\sequence{\sigma}{0} \coloneqq 0$, $\sequence{r_0}{0} \coloneqq 0$, $k \coloneqq 0$, $K \coloneqq \emptyset$, $\sequence{J}{0} \coloneqq \{0\}$.
		\While{true}
		\State
    \zcsetup{reftype=line}
    \label{line:subprob}
		Compute a solution $\sequence{\lambda}{k} \in \bbR^{\abs{\sequence{J}{k}}}$ of
		\begin{align*}
			\argmin_{\lambda \in \bbR^{\abs{\sequence{J}{k}}}}
			\quad
			&
			\frac{1}{2}
			\riemanniannorm[Big]{
				\sum_{j \in \sequence{J}{k}}
				\lambda_j
				\,
				\parallelTransport{\sequence{q}{j}}{\sequence{p}{k}}
				X_{\sequence{q}{j}}
			}^2
			+
			\sum_{j \in \sequence{J}{k}}
			\lambda_j
			\,
			\sequence{e_j}{k}
			+
			\sum_{j \in \sequence{J}{k}}
			\sequence{\lambda_j}{k}
			\,
			\sequence{r_j}{k}
			\\
			\text{\suto}
			\quad
			&
			\sum_{j \in \sequence{J}{k}} \lambda_j
			=
			1
			\quad
			\text{and }
			\lambda_j
			\ge
			0
			\quad \text{for all }
			j \in \sequence{J}{k}
			.
		\end{align*}
		\State
		\label{line:search-direction}
    \zcsetup{reftype=line}
		Set $\sequence{g}{k} \coloneqq \sum_{j \in \sequence{J}{k}} \sequence{\lambda_j}{k} \, \parallelTransport{\sequence{q}{j}}{\sequence{p}{k}} X_{\sequence{q}{j}}$
		\State
		\label{line:epsilon-k}
    \zcsetup{reftype=line}
		Set $\sequence{\varepsilon}{k} \coloneqq \sum_{j \in \sequence{J}{k}} \sequence{\lambda_j}{k} \, \sequence{e_j}{k}$
		\State
		\label{line:sigma-k}
    \zcsetup{reftype=line}
		Set $\sequence{\sigma}{k} \coloneqq \sum_{j \in \sequence{J}{k}} \sequence{\lambda_j}{k} \, \sequence{r_j}{k}$
		\State
		Set $\sequence{d}{k} \coloneqq -\sequence{g}{k}$
		\label{line:search-direction-1}
    \zcsetup{reftype=line}
		\Comment{search direction}
		\State
		Set $\sequence{\xi}{k} \coloneqq - \riemanniannorm{\sequence{g}{k}}^2 - \sequence{\varepsilon}{k} - \sequence{\sigma}{k}$
		\If{$- \sequence{\xi}{k} \le \tol$}
		\label{line:stopping-criterion}
    \zcsetup{reftype=line}
		\Comment{stopping criterion}
		\State \Return $\sequence{p}{k}$
		\EndIf
		\State
		Set $\sequence{t}{k} \coloneqq 1$
		\label{line:initial-step-size}
    \zcsetup{reftype=line}
		\While{$\sequence{q}{k+1} \coloneqq \exponential{\sequence{p}{k}}(\sequence{t}{k} \sequence{d}{k}) \not \in \interior \dom f$, \\
        \hspace*{1.5cm} or $\dist(\sequence{p}{k},\sequence{q}{k+1}) < \sequence{t}{k} \riemanniannorm{\sequence{d}{k}}$}
		\label{line:backtracking-loop}
    \zcsetup{reftype=line}
		\Comment{trial step}
		\State
		Set $\sequence{t}{k} \coloneqq \beta \, \sequence{t}{k}$
		\Comment{backtracking}
		\EndWhile
		\label{line:end-backtracking-loop}
    \zcsetup{reftype=line}
		\State
		Evaluate $X_{\sequence{q}{k+1}} \in \partial f(\sequence{q}{k+1})$
		\If {$f(\sequence{q}{k+1}) \le f(\sequence{p}{k}) + m \, \sequence{t}{k} \, \sequence{\xi}{k}$}
		\label{line:descent-test}
    \zcsetup{reftype=line}
		\State
		\label{line:serious-step}
    \zcsetup{reftype=line}
		Set $\sequence{\bar t}{k} \coloneqq \sequence{t}{k}$ and $\sequence{p}{k+1} \coloneqq \sequence{q}{k+1}$ and $K \coloneqq K \cup \{k\}$
		\Comment{serious step}
		\Else
		\State
		Set $\sequence{\bar t}{k} \coloneqq 0$ and $\sequence{p}{k+1} \coloneqq \sequence{p}{k}$
		\While{
			$
			m
			\,
			\sequence{t}{k}
			\,
			\sequence{\xi}{k}
			\ge
			\riemannian[big]{
				\parallelTransport{\sequence{q}{k+1}}{\sequence{p}{k}}
				X_{\sequence{q}{k+1}}
			}{
				\sequence{t}{k} \,
				\sequence{d}{k}
			}
			-
			\sequence{e_{k+1}}{k}
			-
			\sequence{r_{k+1}}{k}
			$
		}
    \label{line:start-backtracking-loop-null-step}
    \zcsetup{reftype=line}
		\State
		Set $\sequence{t}{k} \coloneqq \beta \, \sequence{t}{k}$
		and $\sequence{q}{k+1} \coloneqq \exponential{\sequence{p}{k}}(\sequence{t}{k} \sequence{d}{k})$
		\EndWhile
		\label{line:end-backtracking-loop-null-step}
    \zcsetup{reftype=line}
		\Comment{null step}
		\EndIf
		\State \label{line:positive-coefficients}
    \zcsetup{reftype=line}
		Set $\sequence{J_+}{k} \coloneqq \setDef[big]{j \in \sequence{J}{k}}{\sequence{\lambda_j}{k} > 0}$ and $\sequence{J}{k+1} \coloneqq \sequence{J_+}{k} \cup \{k+1\}$
		\State
		\label{line:lin-errors}
    \zcsetup{reftype=line}
		Set $\sequence{e_j}{k+1} \coloneqq f(\sequence{p}{k+1}) - f(\sequence{q}{j}) - \riemannian[big]{X_{\sequence{q}{j}}}{\logarithm{\sequence{q}{j}} \sequence{p}{k+1}}$ for all $j~\in~\sequence{J}{k+1}$
		\State
		\label{line:curvature-remainder}
    \zcsetup{reftype=line}
		Set $\sequence{r_j}{k+1} \coloneqq \varrho \, \riemanniannorm[big]{X_{\sequence{q}{j}}} \, \riemanniannorm[big]{\logarithm{\sequence{q}{j}} \sequence{p}{k+1}}$ for all $j~\in~\sequence{J}{k+1}$
		\State
		Set $k \coloneqq k+1$
		\EndWhile
	\end{algorithmic}
\end{algorithm}

We would now like to point out some more details concerning the proposed method.
\begin{enumerate}
	\item
		The backtracking loop in \zcref{line:backtracking-loop} to \zcref{line:end-backtracking-loop} of \zcref{algorithm:RCBM} serves a dual purpose.
		On the one hand, it ensures that the trial and serious iterates remain inside the interior of the domain of~$f$.
		On the other hand, on a compact manifold, it may happen that a geodesic \enquote{wraps around}, so that $\exponential{\sequence{p}{k}}(\sequence{t}{k} \sequence{d}{k})$ is closer to $\sequence{p}{k}$ than $\sequence{t}{k} \riemanniannorm{\sequence{d}{k}}$.
    This situation is avoided with a strategy inspired by \cite[Definition~2.8]{BergmannHerzogSilvaLouzeiroTenbrinckVidalNunez:2021:1}.
		To see that the loop in \zcref{line:backtracking-loop} of \zcref{algorithm:RCBM} is finite, we argue as follows.
		Given that $\sequence{p}{k} \in \interior(\dom f)$ for each $k$, and that $\interior(\dom f)$ is an open set, there exists a ball $\openBall{r}{\sequence{p}{k}} \subseteq \interior(\dom f)$ for some $r > 0$.
		Hence, as long as the minimizer does not lie on the boundary of $\dom f$, the loop stops as soon as $\sequence{t}{k} \riemanniannorm{\sequence{d}{k}} < r$.
		This is guaranteed to happen as the backtracking factor~$\beta$ decreases the step size~$\sequence{t}{k}$.

	\item \label{remark:null-step-condition}
    The loop from \zcref{line:start-backtracking-loop-null-step} to \zcref{line:end-backtracking-loop-null-step} is inspired by \cite[Inequality~(1.30), p.15]{BagirovGaudiosoKarmitsaMakelaTaheri:2020:1}.
		This remark is dedicated to give further motivation for introducing such loop.
		Let $\sequence{q}{k+1}$ be a null step, \ie
		\begin{equation*}
			f(\sequence{q}{k+1})
			>
			f(\sequence{p}{k})
			+
			m
			\,
			\sequence{t}{k}
			\,
			\sequence{\xi}{k}
			.
		\end{equation*}
		Then, by the definition of linearization errors one obtains that
		\begin{align*}
			\sequence{e_{k+1}}{k}
			&
			=
			f(\sequence{p}{k}) - f(\sequence{q}{k+1})
			-
			\riemannian[big]{X_{\sequence{q}{k+1}}}{\logarithm{\sequence{q}{k+1}} \sequence{p}{k}}
			\\
			&
			=
			f(\sequence{p}{k}) - f(\sequence{q}{k+1})
			+
			\riemannian[big]{
				\parallelTransport{\sequence{q}{k+1}}{\sequence{p}{k}}
				X_{\sequence{q}{k+1}}
			}{
				\logarithm{\sequence{p}{k}} \sequence{q}{k+1}
			}
			.
		\end{align*}
		Thus,
		\begin{equation}
			\begin{aligned}
				f(\sequence{q}{k+1}) - f(\sequence{p}{k})
				&
				=
				\riemannian[big]{
					\parallelTransport{\sequence{q}{k+1}}{\sequence{p}{k}}
					X_{\sequence{q}{k+1}}
				}{
					\logarithm{\sequence{p}{k}} \sequence{q}{k+1}
				}
				-
				\sequence{e_{k+1}}{k}
				\\
				&
				=
				-
				\riemannian[big]{
					\parallelTransport{\sequence{q}{k+1}}{\sequence{p}{k}}
					X_{\sequence{q}{k+1}}
				}{
					\sequence{t}{k}
				\sequence{g}{k}}
				-
				\sequence{e_{k+1}}{k}
				,
			\end{aligned}
		\end{equation}
		because $\sequence{q}{k+1} = \exponential{\sequence{p}{k}} (\sequence{t}{k} \sequence{d}{k})$ with $\sequence{d}{k} = - \sequence{g}{k}$ and $\sequence{t}{k} \in (0, 1]$.
		Notice that the second equality above is ensured by the second condition in the backtracking loop; see \zcref{line:backtracking-loop} in \zcref{algorithm:RCBM}.
		Combining the two previous expressions, we get the following strict inequality
		\begin{equation*}
			\riemannian[big]{
				\parallelTransport{\sequence{q}{k+1}}{\sequence{p}{k}}
				X_{\sequence{q}{k+1}}
			}{\sequence{t}{k} \, \sequence{g}{k}}
			<
			-
			m
			\,
			\sequence{t}{k}
			\,
			\sequence{\xi}{k}
			-
			\sequence{e_{k+1}}{k}
			.
		\end{equation*}
		Hence, for null steps, there exists a slack $\alpha \ge 0$ such that
		\begin{equation*}
			\riemannian[big]{
				\parallelTransport{\sequence{q}{k+1}}{\sequence{p}{k}}
				X_{\sequence{q}{k+1}}
			}{\sequence{t}{k} \, \sequence{g}{k}}
			=
			-
			m
			\,
			\sequence{t}{k}
			\,
			\sequence{\xi}{k}
			-
			\sequence{e_{k+1}}{k}
			-
			\alpha
			.
		\end{equation*}
		Thus, the backtracking loop in question aims to find a small enough $\sequence{\bar t}{k} \le \sequence{t}{k}$ such that
		\begin{equation}
			\riemannian[big]{
				\parallelTransport{\sequence{\bar q}{k+1}}{\sequence{p}{k}}
				X_{\sequence{\bar q}{k+1}}
			}{\sequence{\bar t}{k} \, \sequence{g}{k}}
			<
			-
			m
			\,
			\sequence{\bar t}{k}
			\,
			\sequence{\xi}{k}
			-
			\sequence{e_{k+1}}{k}
			-
			\sequence{r_{k+1}}{k}
			,
		\end{equation}
		for $\sequence{\bar q}{k+1} \coloneqq \exponential{\sequence{p}{k}}(-\sequence{\bar t}{k} \sequence{g}{k})$, $X_{\sequence{\bar q}{k+1}} \in \partial f(\sequence{\bar q}{k+1})$, and updated $\sequence{e_{k+1}}{k}$ and $\sequence{r_{k+1}}{k}$ accordingly.
		Observe now that, by the definition of linearization errors, this last inequality implies
		\begin{align*}
			f(\sequence{\bar q}{k+1})
			>
			f(\sequence{p}{k})
			+
			m
			\,
			\sequence{\bar t}{k}
			\,
			\sequence{\xi}{k}
			+
			\sequence{r_{k+1}}{k}
			\ge
			f(\sequence{p}{k})
			+
			m
			\,
			\sequence{\bar t}{k}
			\,
			\sequence{\xi}{k}
			,
		\end{align*}
		which means that the newly found $\sequence{\bar q}{k+1}$ is still a null iterate.
		Finally, we note here that the loop in \zcref{line:end-backtracking-loop-null-step} in \zcref{algorithm:RCBM} has always stopped within a finite number of iterations in our numerical experiments.

	\item \label{remark:generalized-algorithm}
		Following \zcref{remark:retraction-vector-transport-generalization}, the exponential map in~\zcref{line:backtracking-loop} may be replaced by a retraction.
		Analogously, the logarithmic map in \zcref{line:lin-errors} and \zcref{line:positive-coefficients} can be replaced by the (local) inverse retraction.
		Moreover, the parallel transport in \zcref{line:subprob}, and \zcref{line:search-direction} can be replaced by a continuous vector transport~$\vectorTransportSymbol$.
		This amounts to using the quantities $\sequence{\tilde e_j}{k}$ defined in \zcref{eq:inverse-retraction-linearization-errors} in place of $\sequence{e_j}{k}$ in \zcref{line:subprob}, \zcref{line:epsilon-k}, and \zcref{line:lin-errors}.
		Analogously, $\sequence{\tilde r_j}{k}$ defined in \zcref{eq:inverse-retraction-curvature-remainder} is to be used in place of $\sequence{r_j}{k}$ in \zcref{line:subprob}, \zcref{line:sigma-k}, and \zcref{line:curvature-remainder}.
		Additional conditions have to be worked out that ensure $\sequence{\tilde e_j}{k} \ge 0$.
		We can obtain an approximate Riemannian convex bundle method along these lines, but we leave the details to future research.
\end{enumerate}

\section{Convergence of the Method}
\label{section:convergence}

The convergence analysis of \zcref{algorithm:RCBM} follows along the lines of \cite[Chapter~6.7]{GeigerKanzow:2002:1}, and we point out where adjustments are necessary due to the nonzero curvature of the space and the fact that the objective may take the value~$\infty$.
Such adjustments are already necessary in the first step in the analysis, which is an adaptation of~\cite[Lemmas~6.73 and~6.74]{GeigerKanzow:2002:1}.
Here we have to work with the linearization errors $\sequence{e_j}{k}$, the curvature-dependent remainders $\sequence{r_j}{k}$, and with parallel transport.
Throughout, we assume \zcref{assumption:standing-assumptions} to hold.
As is customary, we carry out the convergence analysis with the stopping criterion disabled.

\begin{lemma}
	\label{lemma:algo-numbers-properties}
	The sequences $\sequence[big](){p}{k}$, $\sequence[big](){q}{k}$, $\sequence[big](){g}{k}$, $\sequence[big](){\varepsilon}{k}$, $\sequence[big](){\sigma}{k}$, $\sequence[big](){\xi}{k}$ generated by \zcref{algorithm:RCBM} possess the following properties:
	\begin{enumeratearabic}
		\item \label{item:algo-numbers-properties:1}
			$\sequence{e_j}{k}, \sequence{r_j}{k} \ge 0$, and $\parallelTransport{\sequence{q}{j}}{\sequence{p}{k}} X_{\sequence{q}{j}} \in \partial_{\sequence{e_j}{k} + \sequence{r_j}{k}} f(\sequence{p}{k})$ for all $j \in \sequence{J}{k}$ and $k \in \bbN_0$;

		\item \label{item:algo-numbers-properties:2}
			$\sequence{\varepsilon}{k}, \sequence{\sigma}{k} \ge 0$, and $\sequence{g}{k} \in \partial_{\sequence{\varepsilon}{k} + \sequence{\sigma}{k}} f(\sequence{p}{k})$ for all $k \in \bbN_0$;

		\item \label{item:algo-numbers-properties:3}
			$\sequence{\xi}{k} \le 0$ for all $k \in \bbN_0$;

		\item \label{item:algo-numbers-properties:4}
			$\sequence{p}{k}$ belongs to $\interior \dom f$ for all $k \in \bbN_0$;

		\item \label{item:algo-numbers-properties:5}
			if $\sequence{\xi}{k} = 0$ for some $k \in \bbN_0$, then $\sequence{p}{k}$ is a global minimizer of $f$.
	\end{enumeratearabic}
\end{lemma}
\begin{proof}
	\ref{item:algo-numbers-properties:1}:
	The first statement follows directly by applying the subgradient inequality to $\sequence{e_j}{k}$.
  Furthermore, since $\varrho \ge 0$ holds, it follows that $\sequence{r_j}{k} \ge 0$ for all $j = 0, \ldots, k$ and $k \in \bbN_{0}$.

	For the second statement, consider first the case $k = 0$.
	Then $\sequence{p}{0} = \sequence{q}{0}$ and
	\begin{equation*}
		\sequence{e_0}{0}
		=
		0
		.
	\end{equation*}
	Hence
	\begin{equation*}
		X_{\sequence{p}{0}}
		\in
		\partial f(\sequence{p}{0})
		=
		\partial_0 f(\sequence{p}{0})
		=
		\partial_{\sequence{e_0}{0}} f(\sequence{p}{0})
		.
	\end{equation*}
	For $k \ge 0$ we have
	\begin{equation*}
		\sequence{e_j}{k+1}
		=
		f(\sequence{p}{k+1})
		-
		f(\sequence{q}{j})
		-
		\riemannian[big]{X_{\sequence{q}{j}}}{\logarithm{\sequence{q}{j}} \sequence{p}{k+1}}
		,
	\end{equation*}
	for all $j \in \sequence{J}{k+1}$, where $X_{\sequence{q}{j}} \in \partial f(\sequence{q}{j})$.
	Hence, applying \zcref{ineq:almost-subgradient}, we find
	\begin{equation*}
		f(p)
		\ge
		f(\sequence{p}{k+1})
		+
		\riemannian[big]{\parallelTransport{\sequence{q}{j}}{\sequence{p}{k+1}} X_{\sequence{q}{j}}}{\logarithm{\sequence{p}{k+1}} p}
		-
		\sequence{e_j}{k+1}
		-
		\sequence{r_j}{k+1}
	\end{equation*}
	for all $p \in \dom f$, so that $\parallelTransport{\sequence{q}{j}}{\sequence{p}{k+1}} X_{\sequence{q}{j}} \in \partial_{\sequence{e_j}{k+1} + \sequence{r_j}{k+1}} f(\sequence{p}{k+1})$ for all $j \in \sequence{J}{k+1}$.

	\ref{item:algo-numbers-properties:2}:
	Similarly, by the definitions of $\sequence{\varepsilon}{k}, \sequence{\sigma}{k}$ and $\lambda_j$, we obtain $\sequence{\varepsilon}{k} \ge 0$, and $\sequence{\sigma}{k} \ge 0$ for all $k \in \bbN_0$.
	Now let $p \in \dom f$ and note that, thanks to the previous point,
	\begin{align*}
		f(p)
		&
		=
		\sum_{j \in \sequence{J}{k}} \sequence{\lambda_j}{k}
		\,
		f(p)
		\\
		&
		\ge
		\sum_{j \in \sequence{J}{k}} \sequence{\lambda_j}{k}
		\,
		\paren[big][]{%
			f(\sequence{p}{k})
			+
			\riemannian[big]{
				\parallelTransport{\sequence{q}{j}}{\sequence{p}{k}}
				X_{\sequence{q}{j}}
			}{
				\logarithm{\sequence{p}{k}} p
			}
			-
			\sequence{e_j}{k}
			-
      \sequence{r_j}{k}
		}
		\\
		&
		=
		f(\sequence{p}{k})
		+
		\riemannian[big]{\sequence{g}{k}}{\logarithm{\sequence{p}{k}} p}
		-
		\sequence{\varepsilon}{k}
		-
		\sequence{\sigma}{k}
		.
	\end{align*}
	Therefore $\sequence{g}{k} \in \partial_{\sequence{\varepsilon}{k} + \sequence{\sigma}{k}} f(\sequence{p}{k})$ for all $k \in \bbN_0$.

	\ref{item:algo-numbers-properties:3}:
	Since $\sequence{\xi}{k} = - \riemanniannorm{\sequence{g}{k}}^2 - \sequence{\varepsilon}{k} - \sequence{\sigma}{k}$, we obtain the statement from the previous point.

	\ref{item:algo-numbers-properties:4}:
	This property follows immediately from the condition $\sequence{p}{0} \in \interior \dom f$ and the way subsequent serious iterates are constructed in \zcref{algorithm:RCBM}.

	\ref{item:algo-numbers-properties:5}:
	If $\sequence{\xi}{k} = 0$ for some $k \in \bbN_0$, then $\sequence{g}{k} = 0$, $\sequence{\varepsilon}{k} = 0$, and $\sequence{\sigma}{k} = 0$.
	Therefore
	\begin{align*}
		0
		=
		\sequence{g}{k} \in \partial_0 f(\sequence{p}{k})
		=
		\partial f(\sequence{p}{k})
		.
	\end{align*}
\end{proof}

\begin{remark}
	\label{remark:stopping-criterion-tolerance}
	\zcref{lemma:algo-numbers-properties}~item~\ref{item:algo-numbers-properties:5} gives a condition for stopping \zcref{algorithm:RCBM} motivated by $0 \in \partial f(p^*)$ being a sufficient condition for global optimality, that is also necessary when $p^* \in \interior \dom f$.
	Given a real number $\tol > 0$, the stopping criterion in \zcref{line:stopping-criterion} reads $- \sequence{\xi}{k} \le \tol$.
	This means $\riemanniannorm{\sequence{g}{k}}^2 + \sequence{\varepsilon}{k} + \sequence{\sigma}{k} \le \tol$, implying that $\sequence{p}{k}$ is close to satisfying $0 \in \partial f(\sequence{p}{k})$.
\end{remark}

The following result, which is an adaptation of~\cite[Lemma~6.75]{GeigerKanzow:2002:1}, shows that the convergence of $\sum_{k=0}^{\infty} - \sequence{\bar t}{k} \, \sequence{\xi}{k}$, together with the fact that, if infinitely many serious steps occur, the search directions approach the zero tangent vector and since $\sequence{\varepsilon}{k} + \sequence{\sigma}{k} \to 0$, the $\varepsilon$-subdifferential approaches the subdifferential in the following sense
\begin{equation*}
	\lim_{\varepsilon \to 0} \partial_{\varepsilon} f(p)
	\coloneqq
	\bigcap \setDef[auto]{\partial_{\varepsilon} f(p)}{\varepsilon > 0}
	=
	\partial_0 f(p)
	=
	\partial f(p)
	.
\end{equation*}
\begin{lemma}
	\label{lemma:cauchy-convergence}
	Let $f_\ast \in \bbR$ be such that $f(\sequence{p}{k}) \ge f_\ast$ for all $k \in \bbN_0$.
	Then the sequences $\sequence[big](){p}{k}, \sequence[big](){g}{k}, \sequence[big](){\varepsilon}{k}, \sequence[big](){\sigma}{k}$ generated by \zcref{algorithm:RCBM} are such that
	\begin{enumeratearabic}
		\item \label{item:cauchy-convergence:1}
			$f(\sequence{p}{k})-f(\sequence{p}{k+1}) \to 0$ for $k \to+ \infty$

		\item \label{item:cauchy-convergence:2}
			$\sum_{k=0}^{\infty} \sequence{\bar t}{k} \, \paren[big](){\riemanniannorm{\sequence{g}{k}}^2 + \sequence{\varepsilon}{k} + \sequence{\sigma}{k}} \le \dfrac{f(\sequence{p}{0})-f_\ast}{m}$

		\item \label{item:cauchy-convergence:3}
			if there are infinitely many serious steps $k_n \in K$, then $\sequence{g}{k_n} \to 0$, $\sequence{\varepsilon}{k_n} \to 0$, and $\sequence{\sigma}{k} \to 0$ for $n \to \infty$.
	\end{enumeratearabic}
\end{lemma}
\begin{proof}
  By construction, the sequence $\paren[big](){f(\sequence{p}{k})}$ is monotone non-increasing because $\sequence{\xi}{k} \le 0$.
	Since it is also bounded from below by assumption, it converges, which shows item~\ref{item:cauchy-convergence:1}.

	\ref{item:cauchy-convergence:2}:
	For all $k \in \bbN_0$, hence for serious and null steps, by rearranging the inequality in \zcref{line:descent-test} of \zcref{algorithm:RCBM} we get
	\begin{equation*}
		f(\sequence{p}{k}) - f(\sequence{p}{k+1})
		\ge
		m \, \sequence{\bar t}{k} (- \sequence{\xi}{k})
		,
	\end{equation*}
	since $\sequence{\bar t}{k} \in [0, 1]$ and $\sequence{\bar t}{k} \le \sequence{t}{k}$.
	By summation we obtain
	\begin{equation*}
		f(\sequence{p}{0}) - f(\sequence{p}{k+1})
		\ge
		m \sum_{j=0}^k \sequence{\bar t}{j} \, (- \sequence{\xi}{j})
		.
	\end{equation*}
	Since $f(\sequence{p}{k}) \ge f_\ast$ for all $k \in \bbN_0$, then
	\begin{equation*}
		f(\sequence{p}{0}) - f_\ast
		\ge
		m \sum_{j=0}^k \sequence{\bar t}{j} \, (- \sequence{\xi}{j})
		,
	\end{equation*}
	and by passing to the limit
	\begin{equation*}
		f(\sequence{p}{0}) - f_\ast
		\ge
		m \sum_{j=0}^{\infty} \sequence{\bar t}{j} \, (- \sequence{\xi}{j})
		.
	\end{equation*}
	The claim follows because $\sequence{\xi}{j} = - \riemanniannorm{\sequence{g}{j}}^2 - \sequence{\varepsilon}{j} - \sequence{\sigma}{j}$.

	\ref{item:cauchy-convergence:3}:
	Let $K$ be the set of serious steps defined \zcref{line:serious-step} of \zcref{algorithm:RCBM} and suppose $\abs{K} = + \infty$.
	Then $0 < \sequence{\bar t}{k} \le 1$ if $k \in K$, namely if $k = k_n$ with $k_n < k_{n+1} < \ldots$ for some $n \in \bbN_0$.
	The claim follows from the previous statement, since the series $\sum_{n=1}^{\infty} \sequence{\bar t}{k_n}(- \sequence{\xi}{k_n})$ converges and $\riemanniannorm{\sequence{g}{k_n}}^2 \ge 0$, $\sequence{\varepsilon}{k_n} \ge 0$, and $\sequence{\sigma}{k_n} \ge 0$ by item~\ref{item:algo-numbers-properties:2} of \zcref{lemma:algo-numbers-properties}.
\end{proof}

Next, we show that if \zcref{algorithm:RCBM} generates \emph{infinitely} many serious steps, every accumulation point of the sequence of serious iterates is a minimizer of the objective.
This is an adaptation of~\cite[Lemma~6.76]{GeigerKanzow:2002:1}.

\begin{lemma}
	\label{lemma:infinite-opt-convergence}
	Let $\sequence[big](){p}{k}$ be a sequence generated by \zcref{algorithm:RCBM} and suppose that infinitely many serious steps occur.
	Then every accumulation point $p^* \in \interior(\dom f)$ of $\sequence[big](){p}{k}$ is a minimizer of~$f$.
\end{lemma}
\begin{proof}
	Let $p^* \in \interior(\dom f)$ be an accumulation point for the sequence $\sequence[big](){p}{k}$, with subsequence $\sequence{p}{k_n} \to p^*$.
	Since $f(\sequence{p}{k})$ is monotone non-increasing and there exists a subsequence $f(\sequence{p}{k_n})$ that converges with $\lim_{n \to \infty} f(\sequence{p}{k_n}) \ge f(p^*)$ due to lower semi-continuity, we can conclude $f(\sequence{p}{k}) \ge f(p^*) \eqqcolon f_\ast$ for all $k \in \bbN_0$ and $f(\sequence{p}{k}) \to f_\ast$ for $k \to+ \infty$.
	Since $p^*$ is an accumulation point for $\sequence[big](){p}{k}$ and $\sequence{p}{k}$ does not change through a null step, $p^*$ is an accumulation point for the subsequence $\sequence[big](){p}{k_n}_{k_n \in K}$ as well, with $K$ being the set of indices for the serious steps defined in \zcref{line:serious-step} of \zcref{algorithm:RCBM}.
	By assumption, there exist infinitely many $k_n \in K$ with $k_n < k_{n+1}$ for all $n \in \bbN_0$, so that $\sequence{p}{k_n} \to p^*$ for $n \to \infty$.
	From \zcref{lemma:algo-numbers-properties}~item~\ref{item:algo-numbers-properties:2} we have $\sequence{g}{k} \in \partial_{\sequence{\varepsilon}{k} + \sequence{\sigma}{k}} f(\sequence{p}{k})$ for all $k \in \bbN_0$, so that
	\begin{equation*}
		f(q)
		\ge
		f(\sequence{p}{k_n})
		+
		\riemannian[big]{\sequence{g}{k_n}}{\logarithm{\sequence{p}{k_n}} q}
		-
		\sequence{\varepsilon}{k_n}
		-
		\sequence{\sigma}{k_n}
		\quad
		\text{for all }
		q \in \dom f
		.
	\end{equation*}
	Because $f(\sequence{p}{k}) \ge f_\ast$ holds for all $k \in \bbN_0$, from \zcref{lemma:cauchy-convergence}~item~\ref{item:cauchy-convergence:3} we get $\sequence{g}{k_n} \to 0$, $\sequence{\varepsilon}{k_n} \to 0$, and $\sequence{\sigma}{k_n} \to 0$ for $n \to \infty$.
	By passing to the limit inferior in the previous inequality,
	\begin{align*}
		f(q)
		&
		\ge
		\liminf_{n \to \infty}
		\paren[big][]{%
			f(\sequence{p}{k_n})
			+
			\riemannian[big]{
				\sequence{g}{k_n}
			}{
				\logarithm{\sequence{p}{k_n}} q
			}
			-
			\sequence{\varepsilon}{k_n}
			-
			\sequence{\sigma}{k_n}
		}
		\\
		&
		=
		\liminf_{n \to \infty} f(\sequence{p}{k_n})
		\ge
		f_\ast
	\end{align*}
	for all $q \in \dom f$, where the second inequality is due to the lower semi-continuity of~$f$.
\end{proof}

We now show that if the number of serious iterates generated by \zcref{algorithm:RCBM} is \emph{finite}, then the last such iterate is a minimizer of the objective.
\begin{lemma}
	\label{lemma:finite-opt-convergence}
	Let $\sequence[big](){p}{k}$ be a sequence generated by \zcref{algorithm:RCBM} and suppose that only a finite number of serious steps occur, so that $\sequence{p}{k} = \sequence{p}{k_\ast}$ holds for all $k \ge k_\ast$ for some $k_\ast \in \bbN_0$.
	Then $p^* \coloneqq \sequence{p}{k_\ast}$ is a minimizer of $f$.
\end{lemma}
\begin{proof}
	The main idea in this proof is to bound the objective of the subproblem in \zcref{line:subprob} from above and below using the decomposition of the index set $\sequence{J}{k}$, and to show that this implies the convergence to zero of $\sequence{\xi}{k}$.
	Since $\sequence{p}{k+1} = \sequence{p}{k}$ for all $k \ge k_\ast$, from the update rules of \zcref{algorithm:RCBM} we have that $\sequence{e_j}{k+1} = \sequence{e_j}{k}$, and $\sequence{r_j}{k+1} = \sequence{r_j}{k}$ for all $j \in \sequence{J_+}{k}$ and all $k \ge k_\ast$.
	Let
	\begin{align*}
		\sequence{\eta}{k}
		&
		\coloneqq
		\sum_{j \in \sequence{J_+}{k}}
		\sequence{\lambda_j}{k}
		\,
		\sequence{e_j}{k+1}
		=
		\sequence{\varepsilon}{k}
		,
		\\
		\text{and}
		\quad
		\sequence{\psi}{k}
		&
		\coloneqq
		\sum_{j \in \sequence{J_+}{k}}
		\sequence{\lambda_j}{k}
		\,
		\sequence{r_j}{k+1}
		=
		\sequence{\sigma}{k}
		,
	\end{align*}
	for all $k \ge k_\ast$.
	Let us denote the objective function of the subproblem defined in \zcref{line:subprob} of \zcref{algorithm:RCBM} by
	\begin{equation*}
		\sequence{Q}{k}(\lambda)
		\coloneqq
		\frac{1}{2}
		\riemanniannorm[Big]{
			\sum_{j \in \sequence{J}{k}}
			\lambda_j
			\,
			\parallelTransport{\sequence{q}{j}}{\sequence{p}{k}} X_{\sequence{q}{j}}
		}^2
		+
		\sum_{j \in \sequence{J}{k}}
		\lambda_j
		\,
		\sequence{e_j}{k}
		+
		\sum_{j \in \sequence{J}{k}}
		\lambda_j
		\,
		\sequence{r_j}{k}
		.
	\end{equation*}
	Suppose that $\mu \in [0, 1]$ is an arbitrary number and let $\overline \lambda_j$ be defined for all $j \in \sequence{J}{k} = \sequence{J_+}{k-1} \cup \{k\}$ as
	\begin{equation*}
		\overline \lambda_j
		=
		\begin{cases}
			\mu
			&
			\text{if }
			j = k
			\\
			(1 - \mu)
			\,
			\sequence{\lambda_j}{k-1}
			&
			\text{if }
			j
			\in
			\sequence{J_+}{k-1}
			,
		\end{cases}
	\end{equation*}
	where the numbers $\sequence{\lambda_j}{k-1}$ are the components of the solution $\lambda^{k-1}$ to the subproblem in \zcref{line:subprob} of \zcref{algorithm:RCBM} at the $(k-1)$-st iteration.
	Note that $\overline \lambda_j \ge 0$ for all $j \in \sequence{J}{k}$ and
	\begin{equation*}
		\sum_{j \in \sequence{J}{k}}
		\overline \lambda_j
		=
		\mu
		+
		\sum_{j \in \sequence{J_+}{k-1}}
		(1 - \mu)
		\,
		\sequence{\lambda_j}{k-1}
		=
		\mu
		+
		(1 - \mu)
		\,
		\sum_{j \in \sequence{J}{k-1}}
		\sequence{\lambda_j}{k-1}
		=
		1
	\end{equation*}
	since the numbers $\sequence{\lambda_j}{k-1}$ for $j \in \sequence{J}{k-1} \setminus \sequence{J_+}{k-1}$ are zero for all $k \in \bbN_0$ by definition of $\sequence{J_+}{k}$ and $\sequence{\lambda_j}{k}$.
	This shows that the vector $\overline \lambda \in \bbR^{\abs{\sequence{J}{k}}}$ is admissible for the quadratic program defined in \zcref{line:subprob} of \zcref{algorithm:RCBM}.
	Therefore
	\begin{equation*}
		\sequence{Q}{k}(\sequence{\lambda}{k})
		\le
		\sequence{Q}{k} (\overline \lambda)
		\quad
		\text{for all }
		k \ge k_\ast
		.
	\end{equation*}
	We need to estimate $\sequence{Q}{k}(\overline \lambda)$, and to do this we use the decomposition $\sequence{J}{k} = \sequence{J_+}{k-1} \cup \{k\}$.
	Let $X_{\sequence{q}{j}} \in \partial f(\sequence{q}{j})$ for all $j \in \sequence{J}{k}$ and $k > k_\ast$.
	From the update rules of \zcref{algorithm:RCBM}, and because $\sequence{p}{k} = \sequence{p}{k-1}$ holds for all $k > k_\ast$, we get
	\begin{align*}
		\MoveEqLeft
		\sum_{j \in \sequence{J}{k}}
		\overline \lambda_j
		\,
		\parallelTransport{\sequence{q}{j}}{\sequence{p}{k}} X_{\sequence{q}{j}}
		\\
		&
		=
		\mu
		\,
		\parallelTransport{\sequence{q}{k}}{\sequence{p}{k}} X_{\sequence{q}{k}}
		+
		(1 - \mu)
		\,
		\sum_{j \in \sequence{J}{k-1}}
		\sequence{\lambda_j}{k-1}
		\,
		\parallelTransport{\sequence{q}{j}}{\sequence{p}{k}} X_{\sequence{q}{j}}
		\\
		&
		=
		\mu
		\,
		\parallelTransport{\sequence{q}{k}}{\sequence{p}{k}} X_{\sequence{q}{k}}
		+
		(1 - \mu)
		\,
		\sum_{j \in \sequence{J}{k-1}}
		\sequence{\lambda_j}{k-1}
		\,
		\parallelTransport{\sequence{q}{j}}{\sequence{p}{k-1}} X_{\sequence{q}{j}}
		\\
		&
		=
		\mu
		\,
		\parallelTransport{\sequence{q}{k}}{\sequence{p}{k}} X_{\sequence{q}{k}}
		+
		(1 - \mu)
		\,
		\sequence{g}{k-1}
		,
	\end{align*}
	as well as
	\begin{align*}
		\sum_{j \in \sequence{J}{k}}
		\overline \lambda_j
		\,
		\sequence{e_j}{k}
		&
		=
		\sequence{\overline \lambda}{k}
		\,
		\sequence{e_k}{k}
		+
		\sum_{j \in \sequence{J_+}{k-1}}
		\overline \lambda_j
		\,
		\sequence{e_j}{k}
		\\
		&
		=
		\mu
		\,
		\sequence{e_k}{k}
		+
		\sum_{j \in \sequence{J_+}{k-1}}
		(1 - \mu)
		\,
		\sequence{\lambda_j}{k-1}
		\,
		\sequence{e_j}{k-1}
		\\
		&
		=
		\mu
		\,
		\sequence{e_k}{k}
		+
		(1 - \mu)
		\,
		\sequence{\eta}{k-1}
		.
	\end{align*}
	Analogously, we obtain
	\begin{equation*}
		\sum_{j \in \sequence{J}{k}}
		\overline \lambda_j
		\,
		\sequence{r_j}{k}
		=
		\mu
		\,
		\sequence{r_k}{k}
		+
		(1 - \mu)
		\,
		\sequence{\psi}{k-1}
		.
	\end{equation*}
	By plugging these into the definition of $\sequence{Q}{k}(\overline \lambda)$ we obtain
	\begin{align*}
		\sequence{Q}{k}(\sequence{\lambda}{k})
		&
		\le
		\sequence{Q}{k}(\overline \lambda)
		\\
		&
		=
		\frac{1}{2} \riemanniannorm[big]{
			\mu \, \parallelTransport{\sequence{q}{k}}{\sequence{p}{k}} X_{\sequence{q}{k}}
			+
			(1 - \mu) \, \sequence{g}{k-1}
		}^2
		\\
		&
		\quad
		+
		\mu \, \sequence{e_k}{k}
		+
		(1 - \mu) \, \sequence{\eta}{k-1}
		+
		\mu \, \sequence{r_k}{k}
		+
		(1 - \mu) \, \sequence{\psi}{k-1}
		.
	\end{align*}
	Since $\mu \in [0, 1]$, we see that
	\begin{equation*}
		\sequence{Q}{k}(\sequence{\lambda}{k})
		\le
		\min_{\mu \in [0,1]} \sequence{M}{k}(\mu)
		,
	\end{equation*}
	where
	\begin{align*}
		\sequence{M}{k}(\mu)
		&
		\coloneqq
		\frac{1}{2} \riemanniannorm[Big]{
			\mu \, \parallelTransport{\sequence{q}{k}}{\sequence{p}{k}} X_{\sequence{q}{k}}
			+
			(1 - \mu) \, \sequence{g}{k-1}
		}^2
		\\
		&
		\quad
		+
		\mu \, \sequence{e_k}{k}
		+
		(1 - \mu) \, \sequence{\eta}{k-1}
		+
		\mu \, \sequence{r_k}{k}
		+
		(1 - \mu) \, \sequence{\psi}{k-1}
		.
	\end{align*}
	Let $\sequence{\mu}{k}$ be a solution to the convex, quadratic problem
	\begin{equation*}
		\argmin_{\mu \in [0,1]} \sequence{M}{k}(\mu)
	\end{equation*}
	with optimal value $\sequence{w}k = \sequence{M}{k}(\sequence{\mu}{k})$.
	Then
	\begin{equation}
		\label{ineq:Q-w-M}
		\sequence{w}k
		\le
		\sequence{M}{k}(\mu)
		\quad
		\text{for all }
		\mu \in [0,1]
		.
	\end{equation}
	Hence, for all $k > k_\ast$,
	\begin{align*}
		\sequence{w}k
		&
		\le
		\sequence{M}{k}(0)
		\\
		&
		=
		\frac{1}{2} \riemanniannorm{\sequence{g}{k-1}}^2
		+
		\sequence{\eta}{k-1}
		+
		\sequence{\psi}{k-1}
		\\
		&
		=
		\frac{1}{2} \riemanniannorm[Big]{
			\sum_{j \in \sequence{J}{k-1}}
			\sequence{\lambda_j}{k-1}
			\,
			\parallelTransport{\sequence{q}{j}}{\sequence{p}{k-1}}
			X_{\sequence{q}{j}}
		}^2
		\\
		&
		\quad
		+
		\sum_{j \in \sequence{J}{k-1}}
		\sequence{\lambda_j}{k-1}
		\,
		\sequence{e_j}{k-1}
		+
		\sum_{j \in \sequence{J}{k-1}}
		\sequence{\lambda_j}{k-1}
		\,
		\sequence{r_j}{k-1}
		\\
		&
		=
		\sequence{Q}{k-1}(\lambda^{k-1})
		\\
		&
		\le
		\sequence{w}{k-1}
		.
	\end{align*}
	We can therefore conclude
	\begin{equation}
		\label{ineq:wk-ineq}
		0
		\le
		\sequence{w}k
		\le
		\sequence{w}{k-1}
		\le
		\sequence{w}{k_\ast}
		\quad
		\text{for all }
		k > k_\ast
	\end{equation}
	and
	\begin{equation}
		\label{ineq:gk-norm-ineq}
		\riemanniannorm{\sequence{g}{k}}
		\le
		\sqrt{2 \sequence{w}{k_\ast}}
		,
		\quad
		\sequence{\eta}{k}
		\le
		\sequence{w}{k_\ast}
		,
		\quad
		\text{and }
		\sequence{\psi}{k}
		\le
		\sequence{w}{k_\ast}
		\quad
		\text{for all }
		k \ge k_\ast
		.
	\end{equation}
	Since the $(k-1)$-th step is a null step, from the condition in the loop in \zcref{line:end-backtracking-loop-null-step}, see also \zcref{remark:null-step-condition}, we get
	\begin{align*}
		\riemannian[big]{
			\parallelTransport{\sequence{q}{k}}{\sequence{p}{k}}
			X_{\sequence{q}{k}}
		}{\sequence{t}{k-1} \, \sequence{g}{k-1}}
		&
		<
		m \, \sequence{t}{k-1}
		\,
		\paren[big](){%
			\riemanniannorm{\sequence{g}{k-1}}^2
			+
			\sequence{\eta}{k-1}
			+
			\sequence{\psi}{k-1}
		}
		\\
		&
		\quad
		-
		\sequence{e_k}{k}
		-
		\sequence{r_k}{k}
		,
	\end{align*}
	which implies
	\begin{equation}
		\label{ineq:gk-m-eta}
		\riemannian[big]{
			\parallelTransport{\sequence{q}{k}}{\sequence{p}{k}}
			X_{\sequence{q}{k}}
		}{\sequence{g}{k-1}}
		<
		m
		\,
		\paren[big](){%
			\riemanniannorm{\sequence{g}{k-1}}^2
			+
			\sequence{\eta}{k-1}
			+
			\sequence{\psi}{k-1}
		}
		-
		\sequence{e_k}{k}
		-
		\sequence{r_k}{k}
		.
	\end{equation}

	Let us further estimate $\sequence{M}{k}(\mu)$:
	\begin{align*}
		\sequence{M}{k}(\mu)
		&
		=
		\frac{1}{2} \riemanniannorm[Big]{\mu \, \parallelTransport{\sequence{q}{k}}{\sequence{p}{k-1}} X_{\sequence{q}{k}} + (1 - \mu) \, \sequence{g}{k-1}}^2
		\\
		&
		\quad
		+
		\mu \, \sequence{e_k}{k}
		+
		(1 - \mu) \, \sequence{\eta}{k-1}
		+
		\mu \, \sequence{r_k}{k}
		+
		(1 - \mu) \, \sequence{\psi}{k-1}
		\\
		&
		=
		\frac{1}{2} \riemanniannorm[auto]{\mu \paren[big](){\parallelTransport{\sequence{q}{k}}{\sequence{p}{k-1}} X_{\sequence{q}{k}} - \sequence{g}{k-1} } + \sequence{g}{k-1}}^2
		\\
		&
		\quad
		+
		\mu \, \sequence{e_k}{k}
		+
		(1 - \mu) \, \sequence{\eta}{k-1}
		+
		\mu \, \sequence{r_k}{k}
		+
		(1 - \mu) \, \sequence{\psi}{k-1}
		\\
		&
		=
		\frac{1}{2} \mu^2 \, \riemanniannorm{\parallelTransport{\sequence{q}{k}}{\sequence{p}{k-1}} X_{\sequence{q}{k}} - \sequence{g}{k-1}}^2
		+
		\frac{1}{2} \riemanniannorm{\sequence{g}{k-1}}^2
		\\
		&
		\quad
		+
		\mu \, \paren[auto][]{\riemannian[big]{\parallelTransport{\sequence{q}{k}}{\sequence{p}{k-1}} X_{\sequence{q}{k}}}{\sequence{g}{k-1}} - \riemanniannorm{\sequence{g}{k-1}}^2}
		\\
		&
		\quad
		+
		\mu \, \sequence{e_k}{k} + (1 - \mu) \, \sequence{\eta}{k-1}
		+
		\mu \, \sequence{r_k}{k}
		+
		(1 - \mu) \, \sequence{\psi}{k-1}
		.
	\end{align*}
	From \zcref{ineq:gk-norm-ineq} we see that
	\begin{equation*}
		\frac{1}{2} \riemanniannorm{\sequence{g}{k-1}}^2
		\le
		\frac{1}{2} \riemanniannorm{\sequence{g}{k-1}}^2
		+
		\sequence{\eta}{k-1}
		+
		\sequence{\psi}{k-1}
		\le
		\sequence{w}{k-1}
		,
	\end{equation*}
	and plugging this into the previous inequality and using \zcref{ineq:gk-m-eta} yields
	\begin{align}
		\sequence{M}{k}(\mu)
		&
		\le
		\frac{1}{2} \mu^2 \, \riemanniannorm{\parallelTransport{\sequence{q}{k}}{\sequence{p}{k-1}} X_{\sequence{q}{k}} - \sequence{g}{k-1}}^2
		\notag
		\\
		&
		\quad
		+
		\mu \,
		\paren[Big][]{%
			m
			\,
			\paren[big](){%
				\riemanniannorm{\sequence{g}{k-1}}^2
				+
				\sequence{\eta}{k-1}
				+
				\sequence{\psi}{k-1}
			}
			-
			\sequence{e_k}{k}
			-
			\sequence{r_k}{k}
			-
			\riemanniannorm{\sequence{g}{k-1}}^2
		}
		\notag
		\\
		&
		\quad
		+
		\mu \, \sequence{e_k}{k}
		-
		\mu \, \sequence{\eta}{k-1}
		+
		\mu \, \sequence{r_k}{k}
		-
		\mu \, \sequence{\psi}{k-1}
		+
		\sequence{w}{k-1}
		\notag
		\\
		&
		=
		\frac{1}{2} \mu^2 \riemanniannorm{\parallelTransport{\sequence{q}{k}}{\sequence{p}{k-1}} X_{\sequence{q}{k}} - \sequence{g}{k-1}}^2
		\notag
		\\
		&
		\quad
		-
		\mu \, (1 - m)
		\paren[big](){%
			\riemanniannorm{\sequence{g}{k-1}}^2
			+
			\sequence{\eta}{k-1}
			+
			\sequence{\psi}{k-1}
		}
		+
		\sequence{w}{k-1}
		.
		\label{ineq:final-M-k}
	\end{align}
	The sequences $\sequence[big](){g}{k}$, $\sequence[big](){\eta}{k}$, and $\sequence[big](){\psi}{k}$ are bounded thanks to \zcref{ineq:gk-norm-ineq}.
	Moreover, the subdifferential is locally bounded; see, \eg, \cite[Theorem~4.6]{Udriste:1994:1}.
	Hence, there exists a constant $a > 0$ such that
	\begin{equation}
		\label{ineq:bounded-values}
		\riemanniannorm{X_{\sequence{q}{k}}}
		=
		\riemanniannorm{\parallelTransport{\sequence{q}{k}}{\sequence{p}{k-1}} X_{\sequence{q}{k}}}
		\le
		a
		,
		\quad
		\riemanniannorm{\sequence{g}{k}}
		\le
		a
		,
		\quad
		\sequence{\eta}{k}
		\le
		a
		,
		\quad
		\sequence{\psi}{k}
		\le
		a
		,
	\end{equation}
	for all $k > k_\ast$, so that
	\begin{equation*}
		\riemanniannorm{\parallelTransport{\sequence{q}{k}}{\sequence{p}{k-1}} X_{\sequence{q}{k}} - \sequence{g}{k-1}}^2
		\le
		\paren[big](){\riemanniannorm{\parallelTransport{\sequence{q}{k}}{\sequence{p}{k-1}} X_{\sequence{q}{k}}} + \riemanniannorm{\sequence{g}{k-1}}}^2
		\le
		4 \, a^2
		.
	\end{equation*}
	Plugging this into \zcref{ineq:final-M-k} and defining
	\begin{equation*}
		\sequence{\theta}{k}(\mu)
		\coloneqq
		2 \, a^2 \, \mu^2
		-
		\paren[Big][]{%
			(1 - m)
			\paren[auto](){%
				\riemanniannorm{\sequence{g}{k-1}}^2
				+
				\sequence{\eta}{k-1}
				+
				\sequence{\psi}{k-1}
			}
		}
		\mu
		+
		\sequence{w}{k-1}
		,
	\end{equation*}
	we conclude that
	\begin{equation}
		\label{ineq:M-theta}
		\sequence{M}{k}(\mu)
		\le
		\sequence{\theta}{k}(\mu)
		\quad
		\text{for all }
		k
		>
		k_\ast
		.
	\end{equation}
	Since $\sequence{\theta}{k}$ is a quadratic function of $\mu$, it is straightforward to see that it achieves its minimum value at the point
	\begin{align*}
		\sequence{\mu_\ast}{k}
		&
		=
		\frac{1}{4 \, a^2}
		\paren[Big][]{%
			(1 - m)
			\paren[auto](){%
				\riemanniannorm{\sequence{g}{k-1}}^2
				+
				\sequence{\eta}{k-1}
				+
				\sequence{\psi}{k-1}
			}
		}
		\\
		&
		\le
		\frac{1}{4 \, a^2}
		\paren[auto][]{%
			\riemanniannorm{\sequence{g}{k-1}}^2
			+
			\sequence{\eta}{k-1}
			+
			\sequence{\psi}{k-1}
		}
		\\
		&
		\le
		\frac{a^2 + 2 \, a}{4 \, a^2}
		,
	\end{align*}
	where we used \zcref{ineq:bounded-values}.
	Note that we can choose $a \ge \frac{2}{3}$ without loss of generality, so that $\sequence{\mu_\ast}{k} \le 1$ for all $k > k_\ast$.
	The minimal value of $\sequence{\theta}{k}$ on $[0, 1]$ is then
	\begin{equation}
		\begin{aligned}
			\label{eq:minimal-func-value}
			\sequence{\theta}{k}(\sequence{\mu_\ast}{k})
			&
			=
			2 \, a^2 \, {\sequence{\mu_\ast}{k}}^2
			-
			(1 - m)
			\paren[auto](){%
				\riemanniannorm{\sequence{g}{k-1}}^2
				+
				\sequence{\eta}{k-1}
				+
				\sequence{\psi}{k-1}
			}
			\\
			&
			\quad
			\sequence{\mu_\ast}{k}
			+
			\sequence{w}{k-1}
			.
		\end{aligned}
	\end{equation}
	Hence, since $\sequence{\mu_\ast}{k} \in [0, 1]$, through \zcref{ineq:Q-w-M}, \zcref{ineq:gk-norm-ineq}, \zcref{ineq:M-theta}, and \zcref{eq:minimal-func-value} we get
	\begin{align*}
		\sequence{w}k
		&
		=
		\sequence{M}{k}(\sequence{\mu}{k})
		\le
		\sequence{M}{k}(\sequence{\mu_\ast}{k})
		\le
		\sequence{\theta}{k}(\sequence{\mu_\ast}{k})
		\\
		&
		=
		-
		\frac{1}{8 \, a^2}
		(1 - m)^2
		\paren[auto](){%
			\riemanniannorm{\sequence{g}{k-1}}^2
			+
			\sequence{\eta}{k-1}
			+
			\sequence{\psi}{k-1}
		}^2
		+
		\sequence{w}{k-1}
		,
	\end{align*}
	for all $k > k_\ast$, where the last equality comes from substituting the expression for $\sequence{\mu_\ast}{k}$.
	Hence,
	\begin{equation*}
		\sequence{w}k
		\le
		-
		\frac{1}{8 \, a^2}
		(1 - m)^2
		\paren[auto](){\riemanniannorm{\sequence{g}{k-1}}^2 + \sequence{\eta}{k-1} + \sequence{\psi}{k-1}}^2
		+
		\sequence{w}{k-1}
	\end{equation*}
	for all $k > k_\ast$.
	Rearranging the terms and summing the resulting inequalities for $j = k_\ast+1, \ldots, k+1$ yields
	\begin{equation*}
		C \sum_{j=k_\ast}^k
		\paren[auto](){%
			\riemanniannorm{\sequence{g}{j}}^2
			+
			\sequence{\eta}j
			+
			\sequence{\psi}{j}
		}^2
		\le
		\sequence{w}{k_\ast}
		-
		\sequence{w}{k+1}
		,
	\end{equation*}
	where $C = \dfrac{1}{8a^2} \paren[auto](){1 - m}^2$.
	Now, $\sequence{w}{k+1} \ge 0$ due to \zcref{ineq:wk-ineq}, thus by passing to the limit we have
	\begin{equation*}
		\sum_{j=k_\ast}^{\infty}
		\paren[auto](){%
			\riemanniannorm{\sequence{g}{j}}^2
			+
			\sequence{\eta}{j}
			+
			\sequence{\psi}{j}
		}^2
		<
		+\infty
		.
	\end{equation*}
	Thus $\sequence{g}{k} \to 0$, $\sequence{\eta}{k} = \sequence{\varepsilon}{k} \to 0$, and $\sequence{\psi}{k} = \sequence{\sigma}{k} \to 0$ as $k \to + \infty$, for $k \ge k_\ast$.
	The claim follows by reasoning as in the second part of \zcref{lemma:infinite-opt-convergence}.
\end{proof}

We can now state the following preliminary convergence result.
\begin{theorem}
	\label{theorem:opt-convergence}
	Every accumulation point of the sequence $\sequence[big](){p}{k}$ generated by \zcref{algorithm:RCBM} that lies in $\interior(\dom f)$ is a solution to \zcref{problem:main}.
\end{theorem}
\begin{proof}
	The claim follows by combining \zcref{lemma:infinite-opt-convergence} and \zcref{lemma:finite-opt-convergence}.
\end{proof}
The following lemma guarantees the existence of accumulation points of the sequence of serious iterates of \zcref{algorithm:RCBM} under the assumption that the set of interior solutions of \zcref{problem:main} is nonempty. To show this, we adapted~\cite[Lemma~6.79]{GeigerKanzow:2002:1}, using the function $\zeta_{1, \omega}(\delta)$ defined in \zcref{eq:hessian-eigs}.

\begin{lemma}
	\label{lemma:bounded-min-sequence}
	Let the set $S = \setDef[auto]{p^* \in \interior (\dom f)}{f(p^*) = \inf_{p \in \cM} f(p)}$ of interior solutions to \zcref{problem:main} be nonempty.
	If $p^* \in S$ and $\sequence[big](){p}{k}$ is the sequence of serious iterates generated by \zcref{algorithm:RCBM}, then the following hold:
	\begin{enumeratearabic}
		\item \label{item:bounded-min-sequence:1}
			For all $n \in \bbN_0$ and all $k \ge n$, we have
			\begin{multline*}
				\dist^2(p^*, \sequence{p}{k})
				\le
				\dist^2(p^*, \sequence{p}{n})
				+
				\zeta_{1, \omega}(\delta)
				\sum_{j=n}^k
				\paren[big][]{%
					\dist^2(\sequence{p}{j+1}, \sequence{p}{j})
					+
					2 \, \sequence{\bar t}{j} \,
					\paren[big](){%
						\sequence{\varepsilon}{j}
						+
						\sequence{\sigma}{j}
					}
				}
				.
			\end{multline*}

		\item \label{item:bounded-min-sequence:2}
			$
			\zeta_{1, \omega}(\delta)
			\displaystyle
			\sum_{j=n}^{\infty}
			\paren[big][]{%
				\dist^2(\sequence{p}{j+1}, \sequence{p}{j})
				+
				2 \, \sequence{\bar t}{j} \,
				\paren[big](){%
					\sequence{\varepsilon}{j}
					+
					\sequence{\sigma}{j}
				}
			}
			<
			+ \infty$.

		\item \label{item:bounded-min-sequence:3}
			The sequence $\sequence[big](){p}{k}$ is bounded.
	\end{enumeratearabic}
\end{lemma}
\begin{proof}
	Since $\sequence{g}{k} \in \partial_{\sequence{\varepsilon}{k}+ \sequence{\sigma}{k}} f(\sequence{p}{k})$ for all $k \in \bbN_0$ due to item~\ref{item:algo-numbers-properties:2} of \zcref{lemma:algo-numbers-properties} and since $f(\sequence{p}{k}) \ge f(p^*)$ holds by assumption, then
	\begin{equation*}
		0
		\ge
		f(p^*) - f(\sequence{p}{k})
		\ge
		\riemannian[big]{\sequence{g}{k}}{\logarithm{\sequence{p}{k}} p^*}
		-
		\sequence{\varepsilon}{k}
		-
		\sequence{\sigma}{k}
		,
	\end{equation*}
	which reads
	\begin{equation*}
		\riemannian[big]{\sequence{g}{k}}{\logarithm{\sequence{p}{k}} p^*}
		\le
		\sequence{\varepsilon}{k}
		+
		\sequence{\sigma}{k}
		.
	\end{equation*}
	For $\sequence{\bar t}{k} \in [0, 1]$, by definition,
	\begin{equation*}
		\logarithm{\sequence{p}{k}} \sequence{p}{k+1}
		=
		\sequence{\bar t}{k} \sequence{d}{k} = - \sequence{\bar t}{k} \sequence{g}{k}
	\end{equation*}
	so that
	\begin{equation}
		\label{ineq:tk-epsilonk}
		- \riemannian[big]{\logarithm{\sequence{p}{k}} \sequence{p}{k+1}}{\logarithm{\sequence{p}{k}} p^*}
		\le
		\sequence{\bar t}{k} \,
		\paren[big](){%
			\sequence{\varepsilon}{k}
			+
			\sequence{\sigma}{k}
		}
		.
	\end{equation}
	From \cite[Lemma~6]{ZhangSra:2016:1} we have
	\begin{equation}
		\label{ineq:p*-pk1}
		\begin{aligned}
			\dist^2(p^*, \sequence{p}{k+1})
			&
			\le
			\dist^2(p^*, \sequence{p}{k})
			+
			\zeta_{1, \omega} \paren[big](){\dist(p^*, \sequence{p}{k})}
			\dist^2(\sequence{p}{k}, \sequence{p}{k+1})
			\\
			&
			\quad
			-
			2 \, \riemannian[big]{\logarithm{\sequence{p}{k}} p^*}{\logarithm{\sequence{p}{k}} \sequence{p}{k+1}}
			.
		\end{aligned}
	\end{equation}
	Combining now \zcref{ineq:tk-epsilonk} and \zcref{ineq:p*-pk1} we get
	\begin{multline*}
		\dist^2(p^*, \sequence{p}{k+1})
		\le
		\dist^2(p^*, \sequence{p}{k})
		+
		\zeta_{1, \omega}\paren[big](){\dist(p^*, \sequence{p}{k})}
		\dist^2(\sequence{p}{k+1}, \sequence{p}{k})
		+
		2 \, \sequence{\bar t}{k} \,
		\paren[big](){%
			\sequence{\varepsilon}{k}
			+
			\sequence{\sigma}{k}
		}
		.
	\end{multline*}
	Since $\zeta_{1, \omega}(\dist(p^*, \sequence{p}{k})) \ge 1$ holds, we further obtain
	\begin{multline*}
		\dist^2(p^*, \sequence{p}{k+1})
		\le
		\dist^2(p^*, \sequence{p}{k})
		+
		\zeta_{1, \omega}(\dist(p^*, \sequence{p}{k}))
		\paren[auto][]{%
			\dist^2(\sequence{p}{k+1}, \sequence{p}{k})
			+
			2 \, \sequence{\bar t}{k} \,
			\paren[big](){%
				\sequence{\varepsilon}{k}
				+
				\sequence{\sigma}{k}
			}
		}
		,
	\end{multline*}
	which yields the claim after summation because $\zeta_{1, \omega}(\dist(p^*, \sequence{p}{k})) \le \zeta_{1, \omega}(\delta)$ for all $p^*, \sequence{p}{k} \in \dom f$.

	\ref{item:bounded-min-sequence:2}:
	By definition, for $\sequence{\bar t}{k} \in[0,1]$ we have
	\begin{align*}
		\dist^2(\sequence{p}{k+1},\sequence{p}{k})
		&
		=
		\riemanniannorm{\logarithm[big]{\sequence{p}{k}} \sequence{p}{k+1}}^2
		=
		\paren[big][]{\sequence{\bar t}{k}}^2 \, \riemanniannorm{d^k}^2
		\\
		&
		=
		\paren[big][]{\sequence{\bar t}{k}}^2 \, \riemanniannorm{\sequence{g}{k}}^2
		\le
		2 \paren[big][]{\sequence{\bar t}{k}}^2 \, \riemanniannorm{\sequence{g}{k}}^2
		.
	\end{align*}
	We have
	\begin{align*}
		\zeta_{1, \omega}(\delta)
		\sum_{j=0}^{\infty}
		&
		\paren[auto](){%
			\dist^2(\sequence{p}{j+1},\sequence{p}{j})
			+
			2 \, \sequence{\bar t}{j}
			\paren[big](){%
				\sequence{\varepsilon}{j}
				+
				\sequence{\sigma}{j}
			}
		}
		\\
		&
		\le
		2 \, \zeta_{1, \omega}(\delta)
		\sum_{j=0}^{\infty}
		\paren[Big](){%
			\paren[big][]{\sequence{\bar t}{j}}^2
			\riemanniannorm{\sequence{g}{j}
			}^2
			+
			\sequence{\bar t}{j} \,
			\paren[big](){%
				\sequence{\varepsilon}{j}
				+
				\sequence{\sigma}{j}
			}
		}
		\\
		&
		\le
		2 \, \zeta_{1, \omega}(\delta)
		\sum_{j=0}^{\infty}
		\paren[Big](){%
			\sequence{\bar t}{j}
			\riemanniannorm{\sequence{g}{j}}^2
			+
			\sequence{\bar t}{j} \,
			\paren[big](){%
				\sequence{\varepsilon}{j}
				+
				\sequence{\sigma}{j}
			}
		}
		<
		+ \infty
		,
	\end{align*}
	where the last inequality is due to item~\ref{item:cauchy-convergence:2} of \zcref{lemma:cauchy-convergence}.
	This gives the claim since $\zeta_{1, \omega}(\delta)$ is finite as $\delta < + \infty$ if $\omega < 0$ by \zcref{assumption:standing-assumptions}, and, by \zcref{eq:hessian-eigs}, $\zeta_{1, \omega}(\delta) = 1$ for $\omega \ge 0$ and any $\delta$.

	item~\ref{item:bounded-min-sequence:1} and item~\ref{item:bounded-min-sequence:2} imply the fact that $\sequence[big](){p}{k}$ is bounded.
\end{proof}

We can finally prove the following convergence result, adapted from~\cite[Theorem~6.80]{GeigerKanzow:2002:1}.
\begin{theorem}
	\label{theorem:cbm-convergence}
	Let the set $S = \setDef[auto]{p^* \in \interior (\dom f)}{f(p^*) = \inf_{p \in \cM} f(p)}$ of interior solutions to \zcref{problem:main} be nonempty.
	Every sequence of serious iterates $\sequence[big](){p}{k}$ generated by \zcref{algorithm:RCBM} that does not have accumulation points on the boundary of $\dom f$ converges to a minimizer of~$f$.
\end{theorem}
\begin{proof}
	From item~\ref{item:bounded-min-sequence:3} of \zcref{lemma:bounded-min-sequence}, we have that the sequence $\sequence[big](){p}{k}$ is bounded, hence there exists at least one accumulation point $p^* \in \cM$.
	Since the accumulation points of $\sequence[big](){p}{k}$ do not lie outside of the boundary of $\dom f$ by assumption, we have $p^* \in \interior (\dom f)$, and it is a minimizer of~$f$ because of \zcref{theorem:opt-convergence}, hence $p^* \in S$.

	The idea is to show that the sequence $\sequence[big](){p}{k}$ converges to $p^*$.
	To this end, let $\varepsilon > 0$.
	Since $\sequence[big](){p}{k}$ is bounded, there exists a subsequence that converges to $p^*$, and thus there exists a number $n \in \bbN_0$ such that
	\begin{equation*}
		\dist^2(p^*, \sequence{p}{n})
		\le
		\frac{1}{2} \varepsilon
	\end{equation*}
	and, because of item~\ref{item:bounded-min-sequence:2} of \zcref{lemma:bounded-min-sequence},
	\begin{equation*}
		\zeta_{1, \omega}(\delta)
		\sum_{j=n}^{\infty}
		\paren[Big](){%
			\dist^2(\sequence{p}{j+1}, \sequence{p}{j})
			+
			2 \, \sequence{t}{j} \,
			\paren[big](){%
				\sequence{\varepsilon}{j}
				+
				\sequence{\sigma}{j}
			}
		}
		\le
		\frac{1}{2} \varepsilon
		.
	\end{equation*}
	Substituting these inequalities in the first point of \zcref{lemma:bounded-min-sequence} gives
	\begin{equation*}
		\dist^2(p^*, \sequence{p}{k})
		\le
		\dist^2(p^*, \sequence{p}{n})
		+
		\frac{1}{2} \varepsilon
		\le
		\varepsilon
		,
	\end{equation*}
	for all $k \ge n$.
	The claim follows from $\varepsilon > 0$ being arbitrary.
\end{proof}

\section{Numerical Examples}
\label{section:numerics}

The Riemannian Convex Bundle Method (RCBM, \zcref{algorithm:RCBM}) has been implemented in \julia~1.10 within \manoptjl~0.5.
The implementation\footnote{\href{https://manoptjl.org/stable/solvers/convex_bundle_method/}{manoptjl.org/stable/solvers/convex\_bundle\_method/}} is available with the level of generality provided by retractions and vector transports, as discussed in \zcref{remark:generalized-algorithm} in \zcref{section:algorithm}.
The numerical experiments were performed on a MacBook~Pro~M1, 16~GB~RAM running macOS~14.2.

In the following, we compare the RCBM to two other algorithms for the minimization of non-smooth functions, namely the proximal bundle algorithm (PBA) introduced in~\cite{HoseiniMonjeziNobakhtianPouryayevali:2023:1}, and the subgradient method (SGM) which was introduced for optimization on manifolds in~\cite{FerreiraOliveira:1998:1}.
In one of the experiments, namely the signal denoising problem described in \zcref{subsection:H2-signal}, we also compared the RCBM to the cyclic proximal point algorithm (CPPA) presented in~\cite{Bacak:2014:1}.
All of these methods are available in \manoptjl.
In particular, we provide an implementation of the algorithm described in~\cite{HoseiniMonjeziNobakhtianPouryayevali:2023:1}, which was not publicly available before, to the best of our knowledge.
All algorithms have been used with their default parameters when available, except where explicitly stated otherwise.
As for the PBA, the proximal parameter is initialized at $\mu_0 = \frac{1}{2}$, thereafter being updated with $\mu_i = \log(i+1)$ where $i > 0$ is the $i$-th iteration.
The RCBM, the PBA, and the SGM require a subgradient oracle, whereas the CPPA requires the proximal maps of the objective function.
The tolerance parameter from \zcref{remark:stopping-criterion-tolerance} for stopping the RCBM was set to $\tol = 10^{-8}$, and its maximal number of iterations was set to \num{5000}, in accordance with the default iteration cap of all other algorithms.
For the RCBM, we set the default values for the descent-test parameter to $m = 10^{-3}$, and for the contraction factor $\beta = 0.975$.
For practical reasons, we need to limit the maximal number of elements the bundle can contain.
Hence, we set a maximal bundle size cap of $25$~elements, used for all examples unless specified differently.
When this limit is reached, the oldest element is removed from the bundle unless it coincides with the last serious iterate, in which case it is kept.

Most examples feature $\Omega \le 0$, \ie, the sectional curvature is everywhere non-positive.
This is the case, in particular, on Hadamard manifolds.
These manifolds support real-valued, globally convex functions.
In this case, we may choose the diameter $\delta = \diam(\dom f)$ to be finite but arbitrarily large.
In practice, we discuss the value of $\delta$ directly in the respective experiment subsection.
Two of the examples featured in this section, namely the Riemannian median on the sphere and the spectral Procrustes problem, are meant to investigate the performance of the RCBM outside the setting of geodesic convexity.
This is because the objective functions in these experiments are not necessarily geodesically convex, even when restricted to an arbitrary strongly convex set, although they are still convex in the Euclidean sense.

The runtimes of all experiments in this section have been measured with the \benchmarktools{} package~\cite{ChenRevels:2016:2}.
This package performs either $10^5$ runs of some given code or stops after \SI{5}{\second}, whichever happens first, to then take the median of the measured runtimes.
Since a first run is usually invoking some just-in-time compilation in \julia, we perform one call to the solver before starting the benchmark.
We use this first run also to record the number of iterations and the objective values that are reported in the results.
To solve the quadratic program \zcref{problem:P2}, we use the \ripqp package~\cite{OrbanLeconte:2020:1} with its default parameters, paired with the \quadraticmodels package~\cite{OrbanSiqueira:2019:1} to represent the quadratic problem.

The Riemannian distance plays a key role in our examples.
By defining it as $\dist{q} \colon \cM \to \bbR$ by $\dist{q}{p} = \dist(p, q)$, the subdifferential at $p = q$ is given by
\begin{equation}
	\label{eq:subdifferential-distance}
	\partial \dist{p}{p}
	=
	\setDef[auto]{X_p \in \tangentSpace{p}}{\riemanniannorm{X_p} \le 1},
\end{equation}
see, \eg, \cite[Lemma~4.2]{GrohsHosseini:2015:1}.
We use this fact in computing subgradients for the objective functions in our examples.
In practice, we sample a random unit-length tangent vector $\bar X_p $ at $p$ when $p = q$, so that we get $\bar X_p \in \partial \dist{p}{p}$, and we compute
\begin{equation}
	\label{eq:subgrad-distance}
	\partial \dist{q}{p}
	=
	\paren[big]\{\}{\text{grad} \dist{q}{p}}
	=
	\paren[Big]\{\}{- \frac{\logarithm{p} q}{\dist(p,q)}}
\end{equation}
otherwise.
Further note that, with a similar notation,
\begin{equation}
	\label{eq:subgrad-distance-squared}
	\partial \dist{q}^2(p)
	=
	\paren[big]\{\}{\text{grad} \dist{q}^2(p)}
	=
	\paren[big]\{\}{- 2 \logarithm{p} q}
	.
\end{equation}

\subsection{Riemannian Median}
\label{subsection:riemannian-median}

We discuss three experiments.\footnote{\href{https://juliamanifolds.github.io/ManoptExamples.jl/stable/examples/RCBM-Median/}{juliamanifolds.github.io/ManoptExamples.jl/stable/examples/RCBM-Median/}}
For the first two, let $\cM$ be a Hadamard manifold, \ie, a complete, simply connected Riemannian manifold with non-positive curvature $\Omega \le 0$ everywhere.
Let ${\datasequence{q}{1}, \ldots, \datasequence{q}{N}} \in \cM$ be $N = \num{1000}$~Gaussian random data points, and let $f \colon \cM \to \bbR$ be defined by the globally geodesically convex function
\begin{equation*}
	f(p)
	=
	\sum_{j=1}^N w_j \dist(p, \datasequence{q}{j})
	,
\end{equation*}
where $w_j$, $j = 1, \ldots, N$ are positive weights such that $\sum_{j=1}^N w_j = 1$.
The Riemannian geometric median $p^*$ of the points
\begin{equation*}
	\setDef[big]{\datasequence{q}{1}, \ldots, \datasequence{q}{N}}{\datasequence{q}{j} \in \cM \text{ for all } j = 1, \ldots, N}
\end{equation*}
is then defined as
\begin{equation*}
	p^*
	\coloneqq
	\argmin_{p \in \cM}
	f(p)
	,
\end{equation*}
where equality is justified since $p^*$ is uniquely determined on Hadamard manifolds; see, \eg,~\cite[p.10]{Bacak:2014:1}.
In our experiments, we choose weights $w_j = \frac{1}{N}$.
By \zcref{eq:subgrad-distance}, the subgradients of the objective function are given by
\begin{equation*}
	\partial f(p)
	=
	\frac{1}{N} \sum_{j=1}^{N} \partial \dist{\datasequence{q}{j}}(p)
	,
\end{equation*}
where we choose a random subgradient from the subdifferential.
The initial point for the algorithms is chosen as one of the two points that realize the maximal distance within the dataset and the diameter~$\delta$ is set to two times the maximal distance between the points in the dataset.
The first column of \zcref{table:RM-Hn-comparisons}, \zcref{table:RM-Pn-comparisons}, and \zcref{table:RM-Sn-comparisons} shows the manifold dimension for each run.
Here and in all subsequent tables, we highlight the best results in \textcolor{table-highlight-best}{this color}.

The first experiment is set on the hyperbolic spaces $\cM = \cH^n$, $n = 2^k$ for $k = 1, 2, 5, 10, 15$, equipped with the Minkowski metric of signature $(+, +, \ldots, -)$.
Note that here we have $\omega = \Omega = -1$ everywhere.
The results are summarized in \zcref{table:RM-Hn-comparisons}.
The RCBM and the SGM are largely comparable in terms of their iteration and runtime performance, with the RCBM being faster.
The SGM provides the lowest objective in all runs except for the highest-dimensional experiment.
The optimal objective values attained by the RCBM and the PBA are the same.

\begin{table}[tbp]
  \centering
  \setlength{\tabcolsep}{5.5pt}
  \begin{filecontents*}{data/Riemannian_Median/RCBM-Median_Hn-Comparisons.csv}
2,9,0.00523775,1.0519228600811659,251,0.1320110415,1.0519228806144987,18,0.008112541000000001,1.0474762719782
4,8,0.0046998125,1.075163500336237,230,0.1320907915,1.075163513899843,19,0.0095312915,1.0551779217490966
32,15,0.015195792000000001,1.0855946034628323,234,0.18037358350000002,1.0855946101910428,25,0.0208787705,1.0855948915154165
1024,16,0.284984167,1.0970585121566547,234,4.0077101665,1.0970585211696922,23,0.400037792,1.0970585834918787
32768,16,7.3401725,1.0681033868589755,229,91.28025808400001,1.0681033911649698,21,8.8186905,1.0648750461260832
\end{filecontents*}
\pgfplotstabletypeset[
    col sep = comma,
    every head row/.style = {
      before row = \toprule
      \multicolumn{1}{c}{} & \multicolumn{3}{c}{RCBM} & \multicolumn{3}{c}{PBA} & \multicolumn{3}{c}{SGM}\\
      \cmidrule(lr){2-4}\cmidrule(lr){5-7}\cmidrule(lr){8-10},
      after row = \midrule
    },
    every last row/.style = {after row = \bottomrule},
    display columns/0/.style = {
      column name = Dim.,
      string type,
      column type = {S[table-number-alignment = right, table-format = 5, table-alignment-mode = format]},
    },
    display columns/1/.style = {
      column name = Iter.,
      string type,
      column type = {S[table-number-alignment = right, table-format = 5, table-alignment-mode = format]},
    },
    display columns/2/.style = {
      column name = Time (s),
      string type,
      column type = {S[table-number-alignment = right, table-auto-round = true, scientific-notation = engineering, exponent-mode = scientific, table-format = 1.2e1]},
    },
    display columns/3/.style = {
      column name = Obj.,
      string type,
      column type = {S[table-number-alignment = right, table-auto-round = true, table-format = 1.5, table-alignment-mode = format]},
    },
    display columns/4/.style = {
      column name = Iter.,
      string type,
      column type = {S[table-number-alignment = right, table-format = 5, table-alignment-mode = format]},
    },
    display columns/5/.style = {
      column name = Time (s),
      string type,
      column type = {S[table-number-alignment = right, table-auto-round = true, scientific-notation = engineering, exponent-mode = scientific, table-format = 1.2e1]},
    },
    display columns/6/.style = {
      column name = Obj.,
      string type,
      column type = {S[table-number-alignment = right, table-auto-round = true, table-format = 1.5, table-alignment-mode = format]},
    },
    display columns/7/.style = {
      column name = Iter.,
      string type,
      column type = {S[table-number-alignment = right, table-format = 5, table-alignment-mode = format]},
    },
    display columns/8/.style = {
      column name = Time (s),
      string type,
      column type = {S[table-number-alignment = right, table-auto-round = true, scientific-notation = engineering, exponent-mode = scientific, table-format = 1.2e1]},
    },
    display columns/9/.style = {
      column name = Obj.,
      string type,
      column type = {S[table-number-alignment = right, table-auto-round = true, table-format = 1.5, table-alignment-mode = format]},
    },
    multicolumn names,
    every row 0 column 1/.style={postproc cell content/.style={@cell content={\color{table-highlight-best}##1}}},
    every row 1 column 1/.style={postproc cell content/.style={@cell content={\color{table-highlight-best}##1}}},
    every row 2 column 1/.style={postproc cell content/.style={@cell content={\color{table-highlight-best}##1}}},
    every row 3 column 1/.style={postproc cell content/.style={@cell content={\color{table-highlight-best}##1}}},
    every row 4 column 1/.style={postproc cell content/.style={@cell content={\color{table-highlight-best}##1}}},
    every row 0 column 2/.style={postproc cell content/.style={@cell content={\color{table-highlight-best}##1}}},
    every row 1 column 2/.style={postproc cell content/.style={@cell content={\color{table-highlight-best}##1}}},
    every row 2 column 2/.style={postproc cell content/.style={@cell content={\color{table-highlight-best}##1}}},
    every row 3 column 2/.style={postproc cell content/.style={@cell content={\color{table-highlight-best}##1}}},
    every row 4 column 2/.style={postproc cell content/.style={@cell content={\color{table-highlight-best}##1}}},
    every row 2 column 3/.style={postproc cell content/.style={@cell content={\color{table-highlight-among-best}##1}}},
    every row 3 column 3/.style={postproc cell content/.style={@cell content={\color{table-highlight-among-best}##1}}},
    every row 2 column 6/.style={postproc cell content/.style={@cell content={\color{table-highlight-among-best}##1}}},
    every row 3 column 6/.style={postproc cell content/.style={@cell content={\color{table-highlight-among-best}##1}}},
    every row 0 column 9/.style={postproc cell content/.style={@cell content={\color{table-highlight-among-best}##1}}},
    every row 1 column 9/.style={postproc cell content/.style={@cell content={\color{table-highlight-among-best}##1}}},
    every row 2 column 9/.style={postproc cell content/.style={@cell content={\color{table-highlight-among-best}##1}}},
    every row 3 column 9/.style={postproc cell content/.style={@cell content={\color{table-highlight-best}##1}}},
    every row 4 column 9/.style={postproc cell content/.style={@cell content={\color{table-highlight-best}##1}}},
  ]
  {data/Riemannian_Median/RCBM-Median_Hn-Comparisons.csv}
	\caption{Comparison between the RCBM, the PBA, and the SGM on $\cH^n$ with varying dimensions for the Riemannian median example from \zcref{subsection:riemannian-median}.}
	\label{table:RM-Hn-comparisons}
\end{table}

For the second experiment, we consider symmetric positive definite matrices $\cM = \cP(n)$, $n = 2, 5, 10, 15$, endowed with the affine-invariant metric that makes $\cP(n)$ a Hadamard manifold.
Note that here we have $\omega = -\frac{1}{2},$ and $\Omega = 0$ from \cite[Appendix~I, Proposition~I.1]{CriscitielloBoumal:2020:1}.
A comparison between the algorithms is presented in \zcref{table:RM-Pn-comparisons}.
The PBA outperforms the RCBM in dimension~$15$, while the RCBM is faster than the other algorithms in the remaining experiments, both in terms of number of iterations and CPU runtime.
The RCBM lies between the PBA and SGM in terms of number of iterations, while being more comparable to the former.
All algorithms provide identical objective values up to the five decimal digits shown.

\begin{table}[tbp]
  \centering
  \setlength{\tabcolsep}{5.5pt}
  \begin{filecontents*}{data/Riemannian_Median/RCBM-Median_SPD-Comparisons.csv}
3,43,0.3037513545,0.2608462125867309,57,0.441795896,0.2608462120339797,4629,46.546915542,0.26084621310816863
15,49,2.014069375,0.4365362951414896,75,1.7488511670000002,0.43653629676031297,1727,40.487288583,0.43653629728133464
55,15,1.307492667,0.6180587321332766,89,6.154261875,0.6180587347304607,776,53.362785375,0.6180587347286882
120,6,1.203773396,0.7640309287462824,123,15.406367917,0.7640309361802246,438,53.593232625000006,0.7640309323774311
\end{filecontents*}
\pgfplotstabletypeset[
    col sep = comma,
    every head row/.style = {
      before row = \toprule
      \multicolumn{1}{c}{} & \multicolumn{3}{c}{RCBM} & \multicolumn{3}{c}{PBA} & \multicolumn{3}{c}{SGM}\\
      \cmidrule(lr){2-4}\cmidrule(lr){5-7}\cmidrule(lr){8-10},
      after row = \midrule
    },
    every last row/.style = {after row = \bottomrule},
    display columns/0/.style = {
      column name = Dim.,
      string type,
      column type = {S[table-number-alignment = right, table-format = 5, table-alignment-mode = format]},
    },
    display columns/1/.style = {
      column name = Iter.,
      string type,
      column type = {S[table-number-alignment = right, table-format = 5, table-alignment-mode = format]},
    },
    display columns/2/.style = {
      column name = Time (s),
      string type,
      column type = {S[table-number-alignment = right, table-auto-round = true, scientific-notation = engineering, exponent-mode = scientific, table-format = 1.2e1]},
    },
    display columns/3/.style = {
      column name = Obj.,
      string type,
      column type = {S[table-number-alignment = right, table-auto-round = true, table-format = 1.5, table-alignment-mode = format]},
    },
    display columns/4/.style = {
      column name = Iter.,
      string type,
      column type = {S[table-number-alignment = right, table-format = 5, table-alignment-mode = format]},
    },
    display columns/5/.style = {
      column name = Time (s),
      string type,
      column type = {S[table-number-alignment = right, table-auto-round = true, scientific-notation = engineering, exponent-mode = scientific, table-format = 1.2e1]},
    },
    display columns/6/.style = {
      column name = Obj.,
      string type,
      column type = {S[table-number-alignment = right, table-auto-round = true, table-format = 1.5, table-alignment-mode = format]},
    },
    display columns/7/.style = {
      column name = Iter.,
      string type,
      column type = {S[table-number-alignment = right, table-format = 5, table-alignment-mode = format]},
    },
    display columns/8/.style = {
      column name = Time (s),
      string type,
      column type = {S[table-number-alignment = right, table-auto-round = true, scientific-notation = engineering, exponent-mode = scientific, table-format = 1.2e1]},
    },
    display columns/9/.style = {
      column name = Obj.,
      string type,
      column type = {S[table-number-alignment = right, table-auto-round = true, table-format = 1.5, table-alignment-mode = format]},
    },
    multicolumn names,
    every row 0 column 1/.style={postproc cell content/.style={@cell content={\color{table-highlight-best}##1}}},
    every row 1 column 1/.style={postproc cell content/.style={@cell content={\color{table-highlight-best}##1}}},
    every row 2 column 1/.style={postproc cell content/.style={@cell content={\color{table-highlight-best}##1}}},
    every row 3 column 1/.style={postproc cell content/.style={@cell content={\color{table-highlight-best}##1}}},
    every row 0 column 2/.style={postproc cell content/.style={@cell content={\color{table-highlight-best}##1}}},
    every row 1 column 5/.style={postproc cell content/.style={@cell content={\color{table-highlight-best}##1}}},
    every row 2 column 2/.style={postproc cell content/.style={@cell content={\color{table-highlight-best}##1}}},
    every row 3 column 2/.style={postproc cell content/.style={@cell content={\color{table-highlight-best}##1}}},
    every row 0 column 3/.style={postproc cell content/.style={@cell content={\color{table-highlight-among-best}##1}}},
    every row 1 column 3/.style={postproc cell content/.style={@cell content={\color{table-highlight-among-best}##1}}},
    every row 2 column 3/.style={postproc cell content/.style={@cell content={\color{table-highlight-among-best}##1}}},
    every row 3 column 3/.style={postproc cell content/.style={@cell content={\color{table-highlight-among-best}##1}}},
    every row 0 column 6/.style={postproc cell content/.style={@cell content={\color{table-highlight-among-best}##1}}},
    every row 1 column 6/.style={postproc cell content/.style={@cell content={\color{table-highlight-among-best}##1}}},
    every row 2 column 6/.style={postproc cell content/.style={@cell content={\color{table-highlight-among-best}##1}}},
    every row 3 column 6/.style={postproc cell content/.style={@cell content={\color{table-highlight-among-best}##1}}},
    every row 0 column 9/.style={postproc cell content/.style={@cell content={\color{table-highlight-among-best}##1}}},
    every row 1 column 9/.style={postproc cell content/.style={@cell content={\color{table-highlight-among-best}##1}}},
    every row 2 column 9/.style={postproc cell content/.style={@cell content={\color{table-highlight-among-best}##1}}},
    every row 3 column 9/.style={postproc cell content/.style={@cell content={\color{table-highlight-among-best}##1}}},
  ]
  {data/Riemannian_Median/RCBM-Median_SPD-Comparisons.csv}
	\caption{Comparison between the RCBM, the PBA, and the SGM on $\cP(n)$ with varying dimensions for the Riemannian median example from \zcref{subsection:riemannian-median}.}
	\label{table:RM-Pn-comparisons}
\end{table}

For the third experiment, we consider the sphere $\cS^n$ for $n = 2^k$ for $k = 1, 2, 5, 10, 15$, endowed with the metric inherited from the embedding into $\bbR^{n+1}$.
Note that a major difference here is that this manifold has constant positive sectional curvature $\omega = \Omega = 1$.
In this case, we lose the global convexity of the Riemannian distance and thus of the objective.
Minimizers still exist, but they may, in general, be non-unique.
To circumvent these two issues, we restrict the sampling of data points to a ball of radius $\frac{\pi}{6}$ around the north pole.
We restrict the objective~$f$ to a strongly convex ball of diameter $\delta = \frac{\pi}{3}$ centered at the pole, and extend it to $+\infty$ outside.
By \cite[Proposition~4.2~(ii)]{LiYao:2012:1}, the radial distance is geodesically convex on this ball.

A comparison of the performance of the different algorithms are presented in \zcref{table:RM-Sn-comparisons}.
It is clear from \zcref{table:RM-Sn-comparisons} that the RCBM is fastest both in terms of number of iterations and CPU runtime for the case of dimension $n = 2$, and only in terms of CPU runtime for $n = 2^{15}$.
The PBA is faster in the remaining cases.
In any case, the RCBM and the PBA are largely comparable here, while the SGM was outperformed by both in this experiment.
All algorithms provide essentially indistinguishable objective values.

\begin{table}[tbp]
  \centering
  \setlength{\tabcolsep}{5.5pt}
  \begin{filecontents*}{data/Riemannian_Median/RCBM-Median_Sn-Comparisons.csv}
2,43,0.015813875,0.2588983343894214,71,0.018420292,0.25889833540187207,401,0.097033667,0.25889833463001616
4,74,0.023019708500000003,0.25352537725348034,62,0.016808167000000002,0.2535253775719873,5000,1.332577625,0.253525378184503
32,102,0.043011042,0.25988566942115615,64,0.027273875000000003,0.25988567049286515,231,0.0963391875,0.2598856691132157
1024,103,0.890433688,0.2669928157533156,68,0.622527,0.2669928158844631,185,1.9878565000000001,0.2669928161102596
32768,80,27.199798625000003,0.2593015981501524,65,29.150160375000002,0.2593015971760266,157,205.020359708,0.25930159716205126
\end{filecontents*}
\pgfplotstabletypeset[
    col sep = comma,
    every head row/.style = {
      before row = \toprule
      \multicolumn{1}{c}{} & \multicolumn{3}{c}{RCBM} & \multicolumn{3}{c}{PBA} & \multicolumn{3}{c}{SGM}\\
      \cmidrule(lr){2-4}\cmidrule(lr){5-7}\cmidrule(lr){8-10},
      after row = \midrule
    },
    every last row/.style = {after row = \bottomrule},
    display columns/0/.style = {
      column name = Dim.,
      string type,
      column type = {S[table-number-alignment = right, table-format = 5, table-alignment-mode = format]},
    },
    display columns/1/.style = {
      column name = Iter.,
      string type,
      column type = {S[table-number-alignment = right, table-format = 5, table-alignment-mode = format]},
    },
    display columns/2/.style = {
      column name = Time (s),
      string type,
      column type = {S[table-number-alignment = right, table-auto-round = true, scientific-notation = engineering, exponent-mode = scientific, table-format = 1.2e1]},
    },
    display columns/3/.style = {
      column name = Obj.,
      string type,
      column type = {S[table-number-alignment = right, table-auto-round = true, table-format = 1.5, table-alignment-mode = format]},
    },
    display columns/4/.style = {
      column name = Iter.,
      string type,
      column type = {S[table-number-alignment = right, table-format = 5, table-alignment-mode = format]},
    },
    display columns/5/.style = {
      column name = Time (s),
      string type,
      column type = {S[table-number-alignment = right, table-auto-round = true, scientific-notation = engineering, exponent-mode = scientific, table-format = 1.2e1]},
    },
    display columns/6/.style = {
      column name = Obj.,
      string type,
      column type = {S[table-number-alignment = right, table-auto-round = true, table-format = 1.5, table-alignment-mode = format]},
    },
    display columns/7/.style = {
      column name = Iter.,
      string type,
      column type = {S[table-number-alignment = right, table-format = 5, table-alignment-mode = format]},
    },
    display columns/8/.style = {
      column name = Time (s),
      string type,
      column type = {S[table-number-alignment = right, table-auto-round = true, scientific-notation = engineering, exponent-mode = scientific, table-format = 1.2e1]},
    },
    display columns/9/.style = {
      column name = Obj.,
      string type,
      column type = {S[table-number-alignment = right, table-auto-round = true, table-format = 1.5, table-alignment-mode = format]},
    },
    multicolumn names,
    every row 0 column 1/.style={postproc cell content/.style={@cell content={\color{table-highlight-best}##1}}},
    every row 1 column 4/.style={postproc cell content/.style={@cell content={\color{table-highlight-best}##1}}},
    every row 2 column 4/.style={postproc cell content/.style={@cell content={\color{table-highlight-best}##1}}},
    every row 3 column 4/.style={postproc cell content/.style={@cell content={\color{table-highlight-best}##1}}},
    every row 4 column 4/.style={postproc cell content/.style={@cell content={\color{table-highlight-best}##1}}},
    every row 0 column 2/.style={postproc cell content/.style={@cell content={\color{table-highlight-best}##1}}},
    every row 1 column 5/.style={postproc cell content/.style={@cell content={\color{table-highlight-best}##1}}},
    every row 2 column 5/.style={postproc cell content/.style={@cell content={\color{table-highlight-best}##1}}},
    every row 3 column 5/.style={postproc cell content/.style={@cell content={\color{table-highlight-best}##1}}},
    every row 4 column 2/.style={postproc cell content/.style={@cell content={\color{table-highlight-best}##1}}},
    every row 0 column 3/.style={postproc cell content/.style={@cell content={\color{table-highlight-among-best}##1}}},
    every row 1 column 3/.style={postproc cell content/.style={@cell content={\color{table-highlight-among-best}##1}}},
    every row 2 column 3/.style={postproc cell content/.style={@cell content={\color{table-highlight-among-best}##1}}},
    every row 3 column 3/.style={postproc cell content/.style={@cell content={\color{table-highlight-among-best}##1}}},
    every row 4 column 3/.style={postproc cell content/.style={@cell content={\color{table-highlight-among-best}##1}}},
    every row 0 column 6/.style={postproc cell content/.style={@cell content={\color{table-highlight-among-best}##1}}},
    every row 1 column 6/.style={postproc cell content/.style={@cell content={\color{table-highlight-among-best}##1}}},
    every row 2 column 6/.style={postproc cell content/.style={@cell content={\color{table-highlight-among-best}##1}}},
    every row 3 column 6/.style={postproc cell content/.style={@cell content={\color{table-highlight-among-best}##1}}},
    every row 4 column 6/.style={postproc cell content/.style={@cell content={\color{table-highlight-among-best}##1}}},
    every row 0 column 9/.style={postproc cell content/.style={@cell content={\color{table-highlight-among-best}##1}}},
		every row 1 column 9/.style={postproc cell content/.style={@cell content={\color{table-highlight-among-best}##1}}},
		every row 2 column 9/.style={postproc cell content/.style={@cell content={\color{table-highlight-among-best}##1}}},
		every row 3 column 9/.style={postproc cell content/.style={@cell content={\color{table-highlight-among-best}##1}}},
		every row 4 column 9/.style={postproc cell content/.style={@cell content={\color{table-highlight-among-best}##1}}},
  ]
  {data/Riemannian_Median/RCBM-Median_Sn-Comparisons.csv}
	\caption{Comparison between the RCBM, the PBA, and the SGM on $\cS^n$ with varying dimensions for the Riemannian median example from \zcref{subsection:riemannian-median}.}
	\label{table:RM-Sn-comparisons}
\end{table}

\subsection{Denoising Through Total Variation}
\label{subsection:H2-signal}

Total variation models on manifolds have been researched and implemented in signal and image denoising and restoration problems for manifold valued data~\cite{LellmannStrekalovskiyKoetterCremers:2013:1,WeinmannDemaretStorath:2014:1}.
Let $\cI = \{1, \ldots, n\}$ for $n \ge 1$ be an index set, let $\cN$ be a Hadamard manifold, and consider a manifold-valued signal $q \colon \cI \to \cN$.
We can interpret $q = (\datasequence{q}{1}, \ldots, \datasequence{q}{n})^\transp$ to be a point on the power manifold $\cM = \cN^n$.
Now let $\alpha > 0$ and define the objective $f_q \colon \cM \to \bbR$ as the sum of a data fidelity term $g$ and a total variation term $TV$ as
\begin{equation}\label{eq:tv-objective}
	f_q (p)
	=
	\frac{1}{n}
	\paren[big](){g(p, q) + \alpha \, \operatorname{TV}(p)}
	,
\end{equation}
where
\begin{equation*}
	g(p, q)
	=
	\frac{1}{2} \sum_{i=1}^n \dist{\cN}(\datasequence{p}{i}, \datasequence{q}{i})^2
	\qquad
	\text{and}
	\qquad
	\operatorname{TV}(p)
	=
	\sum_{i=1}^{n-1} \dist{\cN}(\datasequence{p}{i}, \datasequence{p}{i+1}).
\end{equation*}
The non-smoothness of the model comes from the total variation part, where the distance appears with exponent one.
For this, we use \zcref{eq:subdifferential-distance} to compute the subdifferential of the $\operatorname{TV}$ function
\begin{equation*}
	\partial f_q(p)
	=
	\frac{1}{n} \paren[big](){-\logarithm{p} q + \alpha \, \partial \operatorname{TV}(p)}
	,
\end{equation*}
where
\begin{equation*}
	\partial \operatorname{TV}(p)
	=
	\sum_{i=1}^{n-1} \partial \dist{\datasequence{p}{i+1}} (\datasequence{p}{i})
	.
\end{equation*}
As for the proximal maps of the objective function, we refer to~\cite{FerreiraOliveira:2002:1}.
For two signals $p, q \in \cM$ we define the mean squared error as
\begin{equation*}
	E(p, q)
	\coloneqq
	\frac{\dist{\cM}(p,q)}{n}
	=
	\frac{1}{n} \sqrt{\sum_{i=1}^{n} \dist{\cN}^2(\datasequence{p}{i},\datasequence{q}{i})}
	.
\end{equation*}

In this example\footnote{\href{https://juliamanifolds.github.io/ManoptExamples.jl/stable/examples/H2-Signal-TV/}{juliamanifolds.github.io/ManoptExamples.jl/stable/examples/H2-Signal-TV/}}, we consider the task of denoising a measured signal on the hyperbolic manifold $\cN = \cH^2$, a Hadamard manifold.
To this end, let $T > 0$ be a real number, let $\cJ = [a,b] \subseteq \bbR$ be an interval and let $p(t) \in \bbR^2$ be (the graph of) a square wave defined by
\begin{equation*}
	p(t)
	=
	\begin{bmatrix}
		t
		\\
		\sgn(\sin(\frac{2\pi}{T} t))
	\end{bmatrix}
	\quad
	\text{ for all }
	t\in \cI.
\end{equation*}
Let $N > 0$ be the number of points used for the discretization of the interval $\cI$ and let $\tau = \frac{b-a}{N}$ be a step size, so that $\cJ = \bigcup_{i=1}^{N-1} [t_i, t_{i+1}]$ where $t_1 = a$, $t_N = b$ and $t_{i+1} = t_i + \tau$.
We thus have a family of points $\datasequence{p}{i} = p(t_i)_{i=1}^N$ and we select
\begin{equation*}
	\cS
	=
	\paren[big]\{\}{\datasequence{p}{1}}
	\cup
	\setDef[big]{\datasequence{p}{i}, \datasequence{p}{i+1}}{\datasequence{p_2}{i} \neq \datasequence{p_2}{i+1} \text{ for } i = 2, \ldots, N-1}
	\cup
	\paren[big]\{\}{\datasequence{p}{N}}
	,
\end{equation*}
namely the first and last point, together with the points immediately next to the jumps.
To construct points on $\cH^2$ we use the isometry $\varphi \colon \bbR^2 \to \cH^2$ given by
\begin{equation*}
	\varphi(p)
	=
	\begin{bmatrix}
		p_1 \\
		p_2 \\
		\sqrt{\riemanniannorm{p}^2 + 1}
	\end{bmatrix}.
\end{equation*}
We order the points from $\cS$ by time~$t$ into a sequence $(\datasequence{s}{1},\ldots,\datasequence{s}{M})$, where $\datasequence{s}{j} = \varphi(\datasequence{p}{j})$ for the corresponding $\datasequence{p}{j} \in \cS$. Note that by construction, $M$ is even.

Finally, we sample $u = \floor[auto]{\frac{NT}{2(b - a)}}$ points from each shortest geodesic connecting the points $\datasequence{s}{2i-1}$ and $\datasequence{s}{2i}$ for $i = 1, 2, \ldots, \frac{M}{2}$.
In our example, we set $a = -6$, $b = 6$, \ie, $\cJ = [-6,6]$, $T = 3$, and we take $N = 496$~points, which yields a signal $q \in \cM$ of length $n = \frac{uM}{2} = 496$.
We can thus define the power manifold
\begin{equation*}
	\cM
	=
	(\cH^2)^{496}
	.
\end{equation*}
The manufactured signal is a point $q \in \cM$ as shown in \zcref{figure:noise}.
We then generate a noisy signal $\bar q \in\cM$ by adding Gaussian noise to the signal, \ie, we set $\datasequence{\bar q}{i} = \exponential{\datasequence{\bar q}{i}} X_{\datasequence{\bar q}{i}}$ for $i = 1, \ldots, n$, where $X \in \tangentSpace{q}$ has a standard deviation of $\sigma = 0.1$ as in \zcref{figure:noise}.
This noisy signal is used in the data fidelity term $g(p, \bar q)$, and it also serves as the starting point for the algorithms.
The diameter~$\delta$ is set to three times the maximal distance between each of the noisy points in $\cH^2$.
We fix the TV parameter to be $\alpha = 0.5$ for this experiment.

We notice that in this experiment all algorithms reached their iteration cap.
The minimizer $q^* \in \cM$ obtained from minimizing the objective function in \zcref{eq:tv-objective} with the RCBM is shown in \zcref{figure:optimum}.
The other algorithms yield results that are visually very close to the solution obtained by the RCBM.

\begin{figure}
	\pgfplotstableread[col sep = comma]{data/H2_Square_Wave/H2-Signal-TV-data.csv}\datatable
	\pgfplotstableread[col sep = comma]{data/H2_Square_Wave/H2-Signal-TV-geodesics.csv}\datatableGeo
	\pgfplotstableread[col sep = comma]{data/H2_Square_Wave/H2-Signal-TV-noise.csv}\datatableNoise
	\pgfplotstableread[col sep = comma]{data/H2_Square_Wave/H2-Signal-TV-bundle_optimum.csv}\datatableRecon
	\begin{subfigure}[T]{0.48\textwidth}
		\centering
		\begin{tikzpicture}
			\begin{axis}[
				axis equal,
				width = 1.1\textwidth,
				xmin = -1.01, xmax = 1.01,
				ymin = -1.01, ymax = 1.01,
				axis lines = none,
				]
				\addplot[
					TolVibrantGray,
					thin,
					dashed
					] table[x = x, y = y] {\datatableGeo};
				\addplot[
					only marks,
					mark options = {
						draw opacity = 0,
						fill = gray,  
						fill opacity = 0.5,
						mark size = 0.5pt,  
					},
					] table[x = x, y = y] {\datatable};
				\addplot[
					only marks,
					mark options = {
						draw opacity = 0,
						fill = TolVibrantTeal,  
						mark size = 0.5pt,  
					},
					] table[x = x, y = y] {\datatableNoise};

				\draw[gray] (axis cs: 0, 0) circle [radius = 100];
			\end{axis}
		\end{tikzpicture}
		\caption{Artificial signal $q \in (\cH^2)^{496}$ in gray, and noisy data $\bar q \in (\cH^2)^{496}$ in teal.}
		\label{figure:noise}
	\end{subfigure}%
	\hfill
	\begin{subfigure}[T]{0.48\textwidth}
		\centering
		\begin{tikzpicture}
			\begin{axis}[
				axis equal,  
				width = 1.1\textwidth,  
				xmin = -1.01, xmax = 1.01,
				ymin = -1.01, ymax = 1.01,
				axis lines = none, 
				]
				\addplot[
					TolVibrantGray,
					thin,
					dashed
					] table[x = x, y = y] {\datatableGeo};
				\addplot[
					only marks,
					mark options = {
						draw opacity = 0,
						fill = gray,  
						fill opacity = 0.5,
						mark size = 0.5pt,  
					},
					] table[x = x, y = y] {\datatable};
				\addplot[
					only marks,
					mark = *,  
					mark options = {
						draw opacity = 0,
						fill = TolBrightCyan,  
						mark size = 0.5pt,  
					},
					] table[x = x, y = y] {\datatableRecon};
				\draw[gray] (axis cs: 0, 0) circle [radius = 100];
			\end{axis}
		\end{tikzpicture}
		\caption{Denoised reconstruction $q^* \in (\cH^2)^{496}$ in cyan.}
		\label{figure:optimum}
	\end{subfigure}
	\caption{Denoising a signal on $(\cH^2)^{496}$ using the RCBM, \zcref{algorithm:RCBM}. $\operatorname{TV}$ regularization example from \zcref{subsection:H2-signal}.}
\end{figure}

\zcref{table:H2-comparisons} shows how the different algorithms compare to one another in terms of numbers of iterations, CPU runtime in seconds, final objective value, and mean squared error.
Evidently, the RCBM is outperformed in terms of CPU runtime.
Furthermore, the RCBM provides the second-lowest objective value, while the CPPA attains the lowest objective value.
The SGM attains the lowest mean error $E(q^*,q)$, followed by the PBA and the RCBM.

\begin{table}[tbp]
	\begin{subtable}{\textwidth}
		\centering
		\begin{filecontents*}{data/H2_Square_Wave/H2-Signal-TV-comparisons.csv}
Algorithm,Iterations,Time (s),Objective,Error
RCBM,5000,13.892251291000001,0.14023177718722205,0.013691922102575364
PBA,5000,9.519102625,0.14288675212437305,0.013032238475490185
SGM,5000,7.896526708000001,0.14621639996960314,0.0124604833028011
CPPA,5000,3.7385719795,0.13191270168001284,0.017361218377745933
\end{filecontents*}
\pgfplotstabletypeset[col sep = comma,
		display columns/0/.style = {string type, column type = {r}},
		display columns/1/.style = {
			column name = Iter.,
			string type,
			column type = {S[table-number-alignment = right, table-format = 5, table-alignment-mode = format]},
		},
		display columns/2/.style = {
			column name = Time (sec.),
			string type,
			column type = {S[table-number-alignment = right, table-auto-round = true, table-format = 2.3, table-alignment-mode = format]},
		},
		display columns/3/.style = {
			column name = Objective,
			string type,
			column type = {S[table-number-alignment = right, table-auto-round = true, scientific-notation = engineering, exponent-mode = scientific, table-format = 1.4e1, table-column-width = 2.2cm]},
		},
		display columns/4/.style = {
			column name = Error,
			string type,
			column type = {S[table-number-alignment = right, table-auto-round = true, scientific-notation = engineering, exponent-mode = scientific, table-format = 1.4e1, table-column-width = 2.2cm]},
		},
		every head row/.style = {
			before row = \toprule,
			after row = \midrule
		},
		every last row/.style = {after row = \bottomrule},
		multicolumn names,
		every row 0 column 1/.style={postproc cell content/.style={@cell content={\color{table-highlight-best}##1}}},
		every row 1 column 1/.style={postproc cell content/.style={@cell content={\color{table-highlight-best}##1}}},
		every row 2 column 1/.style={postproc cell content/.style={@cell content={\color{table-highlight-best}##1}}},
		every row 3 column 1/.style={postproc cell content/.style={@cell content={\color{table-highlight-best}##1}}},
		every row 3 column 2/.style={postproc cell content/.style={@cell content={\color{table-highlight-best}##1}}},
		every row 3 column 3/.style={postproc cell content/.style={@cell content={\color{table-highlight-best}##1}}},
		every row 2 column 4/.style={postproc cell content/.style={@cell content={\color{table-highlight-best}##1}}},
		]{data/H2_Square_Wave/H2-Signal-TV-comparisons.csv}
	\end{subtable}
	\caption{Comparison of the four algorithms on $(\cH^2)^{496}$ for a $\operatorname{TV}$ parameter of $\alpha = 0.5$. $\operatorname{TV}$ regularization example from \zcref{subsection:H2-signal}.}
	\label{table:H2-comparisons}
\end{table}

\subsection{Spectral Procrustes}
\label{subsection:spectral-procrustes}

We consider another example on a positively curved manifold and non-convex objective.
Given two matrices $A, B \in \bbR^{n \times d}$ we aim to solve the spectral Procrustes problem (see, \eg,~\cite{Cape:2020:1})
\begin{equation}
	\label{problem:spectral-procrustes}
	\argmin_{p \in \SpecialOrthogonal{d}} f(p) \coloneqq \norm{A - B \, p}_2
	,
\end{equation}
where $\SpecialOrthogonal{d}$ is equipped with the standard bi-invariant metric, and $\norm{\,\cdot\,}_2$ denotes the spectral norm of a matrix, \ie, its largest singular value.
\zcref{problem:spectral-procrustes} aims to find the rotation matrix $p \in \bbR^{d \times d}$ which best aligns the columns of~$B$ with those of~$A$.
Note that the spectral norm is convex in the Euclidean sense, but not geodesically convex on $\SpecialOrthogonal{d}$.

To obtain subdifferential information of $f$ in the Riemannian case, we project Euclidean convex subgradients of $f$ onto
\begin{equation*}
	\tangentSpace{p}[\SpecialOrthogonal{d}]
	=
	\setDef[auto]{A \in \bbR^{d \times d}}{p A^\transp + A p^\transp = 0, \; \trace(p^{-1}A) = 0}
	.
\end{equation*}
We use
\begin{equation*}
	\proj{p}(- B^\transp U V^\transp)
	\in
	\tangentSpace{p}[\SpecialOrthogonal{d}]
\end{equation*}
as a substitute for $\partial f(p)$, where $U$ and $V$ are left and right singular vectors, respectively, corresponding to the largest singular value of $A - B \, p$.

The sectional curvature of the orthogonal group $\Orthogonal{d}$ at $p \in \Orthogonal{d}$ is given by
\begin{equation*}
	K_p(X_p, Y_p)
	=
	\frac{1}{4}
	{\riemanniannorm{[X_p, Y_p]}}_F^2
\end{equation*}
for any two orthonormal tangent vectors $X_p, Y_p \in \tangentSpace{p}[\Orthogonal{d}]$, where
\begin{equation*}
	\tangentSpace{p}[\Orthogonal{d}]
	=
	\setDef[auto]{A \in \bbR^{d \times d}}{pA^\top + Ap^\top = 0}
	,
\end{equation*}
see, \eg, \cite[Corollary~3.19]{CheegerEbin:2008:1}.
Here $[X_p, Y_p] = X_p \, Y_p - Y_p \, X_p$ denotes the commutator of $X_p$ and $Y_p$.
\cite[Lemma~2.5]{Ge:2014:1} implies that the curvature of the orthogonal group $\Orthogonal{d}$ is bounded below and above by $\omega = 0$ and $\Omega = \frac{1}{4}$ respectively.
Hence, since $\SpecialOrthogonal{d}$ is the connected component of $\Orthogonal{d}$ that contains the identity, we can use the same bounds on $\SpecialOrthogonal{d}$.

For this experiment\footnote{\href{https://juliamanifolds.github.io/ManoptExamples.jl/stable/examples/Spectral-Procrustes/}{juliamanifolds.github.io/ManoptExamples.jl/stable/examples/Spectral-Procrustes/}}, we set $n = \num{1000}$, $d = 250$, and we initialize two random matrices $A, B \in \bbR^{n \times d}$.
As a starting point for the algorithms, we employ the minimizer of the \emph{orthogonal} Procrustes objective, $\norm{A - B \, p}_F$ , where the minimization is with respect to the Frobenius norm.
The solution to this problem is given in closed form as
\begin{equation*}
	\bar p
	=
	U
	V^\transp
	,
\end{equation*}
see, \eg,~\cite{Schoenemann:1966:1}.
We then set $\sequence{p}{0} = \proj{\SpecialOrthogonal{d}}(\bar p)$ as the starting point for all algorithms.
For RCBM, we set a diameter of $\delta = \frac{\pi}{3 \, \sqrt{\Omega}}$.
The results are summarized in \zcref{table:SP-comparisons} and \zcref{figure:SP-plot}.

\begin{table}[tbp]
	\begin{subtable}{\textwidth}
		\centering
		\begin{filecontents*}{data/Spectral_Procrustes/Spectral-Procrustes-comparisons.csv}
Algorithm,Iterations,Time (s),Objective
RCBM,99,102.03617616700001,235.45980372771464
PBA,31,5.804896042,235.4598037279529
SGM,5000,402.738980167,235.4598037301118
\end{filecontents*}
\pgfplotstabletypeset[col sep = comma,
		display columns/0/.style = {string type, column type = {r}},
		display columns/1/.style = {
			column name = Iter.,
			string type,
			column type = {S[table-number-alignment = right, table-format = 5, table-alignment-mode = format]},
		},
		display columns/2/.style = {
			column name = Time (sec.),
			string type,
			column type = {S[table-number-alignment = right, table-auto-round = true, table-format = 3.2]},
		},
		display columns/3/.style = {
			column name = Objective,
			string type,
			column type = {S[table-number-alignment = right, table-auto-round = true, table-format = 3.8]},
		},
		every head row/.style = {
			before row = \toprule,
			after row = \midrule
		},
		every last row/.style = {after row = \bottomrule},
		multicolumn names,
		every row 1 column 1/.style={postproc cell content/.style={@cell content={\color{table-highlight-best}##1}}},
		every row 1 column 2/.style={postproc cell content/.style={@cell content={\color{table-highlight-best}##1}}},
		every row 2 column 3/.style={postproc cell content/.style={@cell content={\color{table-highlight-among-best}##1}}},
		every row 1 column 3/.style={postproc cell content/.style={@cell content={\color{table-highlight-among-best}##1}}},
		every row 0 column 3/.style={postproc cell content/.style={@cell content={\color{table-highlight-among-best}##1}}},
		]{data/Spectral_Procrustes/Spectral-Procrustes-comparisons.csv}
	\end{subtable}
	\caption{Comparison of the three algorithms on $\SpecialOrthogonal{250}$. Spectral Procrustes example from \zcref{subsection:spectral-procrustes}.}
	\label{table:SP-comparisons}
\end{table}

\begin{figure}
	\pgfplotstableread[col sep = comma]{data/Spectral_Procrustes/Spectral-Procrustes_RCBM-result.csv}\rcbmdata
	\pgfplotstableread[col sep = comma]{data/Spectral_Procrustes/Spectral-Procrustes_PBA-result.csv}\pbadata
	\pgfplotstableread[col sep = comma]{data/Spectral_Procrustes/Spectral-Procrustes_SGM-result.csv}\sgmdata
	\centering
	\begin{tikzpicture}[scale = 0.8]
		\begin{axis}[
			width = 0.9\textwidth,  
			height = 0.33\textheight,
			xmode = log, ymode = normal,
			xmin = 0.5, xmax = 5000,
			ymin = 235.3, ymax = 239.5,
			grid = major,
			xtick = {1e0,1e1,1e2,1e3},
			xticklabels = {1,10,100,1000},
			xlabel = {$k$},
			x label style = {at = {(axis description cs:1.0,0.075)},anchor = north},
			y label style = {at = {(axis description cs:0.125,1.0)},rotate = -90,anchor = west},
			ylabel = {$f(p)$},
			axis lines = left,
			]
			\addplot[
				TolVibrantBlue,
				line width = 1.2pt,
				] table[x = i, y = cost] {\sgmdata};
			\addlegendentry{SGM}
			\addplot[
				TolVibrantMagenta,
				line width = 1.2pt,
				] table[x = i, y = cost] {\pbadata};
			\addlegendentry{PBA}
			\addplot[
				TolVibrantOrange,
				line width = 1.3pt,
				] table[x = i, y = cost] {\rcbmdata};
			\addlegendentry{RCBM}
		\end{axis}
	\end{tikzpicture}
	\caption{Objective versus iteration plot for the spectral Procrustes example from \zcref{subsection:spectral-procrustes}.}
	\label{figure:SP-plot}
\end{figure}%

We see that in this experiment, the PBA outperforms the other two algorithms.
The final objective values are all the same.
From \zcref{figure:SP-plot}, we can further evaluate the development of the cost during the iterations of all three algorithms.
While the RCBM requires more time per iteration, its objective profile stays en par with the PBA, while the SGM is heavily affected by a zig-zagging behavior that leads to a higher number of iterations.

\section{Conclusion}
\label{section:conclusion}

We introduced a new bundle method for solving smooth and non-smooth geodesically convex optimization problems on Riemannian manifolds with bounded curvature.
By formulating a model based on the convex hull of previously collected subgradients, we generalized the dual formulation to minimize the objective function.
We proved convergence of the method to a global minimizer under mild conditions.
The performance of this new algorithm, termed the Riemannian Convex Bundle Method (RCBM), was studied and showcased through various numerical examples implemented using the \julia package \manoptjl, by comparing it to other existing non-smooth optimization algorithms such as the proximal bundle algorithm (PBA), the subgradient method (SGM), and the cyclic proximal point algorithm (CPPA).
For the Riemannian median, our RCBM outperforms the SGM and the PBA both in terms of the number of iterations as well as in runtime on Hadamard manifolds, and the same holds on the sphere for the experiments in dimensions $2$ and $2^{15}$.
Finally, the RCBM may yield meaningful results even when the objective function is not geodesically convex, as shown in the spectral Procrustes experiment.
Further research directions include a better strategy to eliminate unnecessary subgradients from the bundle and investigating the necessary diameter parameter in more detail.
Moreover, one further line of research is the investigation of non-asymptotic complexity bounds for the RCBM.
We also aim to apply the algorithm to other structured, non-smooth optimization problems appearing in applications.

%% file: appendix.tex
\section{Hessian Bound}
\label{section:curvature-bound}

We elaborate on how the argument presented in \cite[Appendix~C]{AlimisisOrvietoBecigneulLucchi:2021:1} was adapted to prove \zcref{theorem:eps-approx}.
Using the notation introduced in said theorem, let $\delta = \diam(\dom f) < +\infty$, $p \in \dom f$, $X_p \in \partial f(p)$, and $\sequence{p}{j} \in \dom f$ and $X_{\sequence{p}{j}} \in \partial f(\sequence{p}{j})$ for all $j = 1, \ldots, k$, $k \in \bbN$.
Let $g \colon [0, 1] \to \bbR$ be a function defined by
\begin{equation}
	g(t)
	=
	\riemannian{X_{\sequence{p}{j}}}
	{
		\parallelTransport{\gamma(t)}{\sequence{p}{j}}
		\logarithm{\gamma(t)}p
	}
	,
\end{equation}
where $\geodesic<s>(t) \colon [0, 1] \to \cM$ is the minimal geodesic arc connecting $\geodesic<s>(0) = \sequence{p}{j}$ to $\geodesic<s>(1) = \sequence{p}{k}$.
The rest of the argument applies almost verbatim as follows.
By applying the mean value theorem to $g(t)$, one can write
\begin{equation}\label{eq:inner-product-hessian}
	\riemannian{X_{\sequence{p}{j}}}
	{
		\logarithm{\sequence{p}{j}}p
		-
		\parallelTransport{\sequence{p}{k}}{\sequence{p}{j}}
		\logarithm{\sequence{p}{k}}p
	}
	=
	\riemannian{X_{\sequence{p}{j}}}
	{
		\cH_{t_0}
		\logarithm{\sequence{p}{j}}\sequence{p}{k}
	}
\end{equation}
for some $t_0 \in (0, 1)$, where
\begin{equation*}
	\cH_t
	\coloneqq
	-
	\parallelTransport{\geodesic<s>(t)}{\sequence{p}{j}}
	\Hess{\geodesic<s>(t)}{-\frac{1}{2}\dist^2(\geodesic<s>(t), p)}
	\parallelTransport{\sequence{p}{j}}{\geodesic<s>(t)}
	,
\end{equation*}
is a family of operators $\cH_t \colon \tangentSpace{\sequence{p}{j}} \to \tangentSpace{\sequence{p}{j}}$ for $t \in [0, 1]$, and $\Hess{\geodesic<s>(t)}{f} \colon \tangentSpace{\geodesic<s>(t)} \to \tangentSpace{\geodesic<s>(t)}$ is the Riemannian Hessian of a $\cC^2$ function $f$ at the point ${\geodesic<s>(t)} \in \cM$.
From \cite[Appendix~D]{AlimisisOrvietoBecigneulLucchi:2021:1}
\begin{equation}\label{ineq:minmax-eigs}
	\norm{\cH_t - \id}
	\le
	\max
	\{
		\zeta_{1, \omega}(\dist(\geodesic<s>(t_0), \, p)) - 1
		,
		1 - \zeta_{2, \Omega}(\dist(\geodesic<s>(t_0), \, p))
	\}
	,
\end{equation}
where the functions $\zeta_{1, \omega}$ and $\zeta_{2, \Omega}$ are, respectively, upper and lower bounds to the largest and smallest eigenvalues of
\begin{equation*}
	-\Hess{\geodesic<s>(t)}{
		-\frac{1}{2}
		\dist^2(
			\geodesic<s>(t)
			,
			\,
			p
		)
	}
	,
\end{equation*}
and they are defined in \zcref{eq:hessian-eigs}.
This comes from $\parallelTransport{\sequence{p}{j}}{\geodesic<s>(t)} = (\parallelTransport{\geodesic<s>(t)}{\sequence{p}{j}})^{-1}$.
Let now $r \in \bbR$ and consider $\varrho_{\omega, \Omega}(r)$ as defined in \zcref{eq:rho-definition}.
For all $p \in \dom f$ we have
\begin{align*}
	- \riemannian{X_{\sequence{p}{j}}}
	{
		\logarithm{\sequence{p}{j}}p
		&
		-
		\parallelTransport{\sequence{p}{k}}{\sequence{p}{j}}
		\logarithm{\sequence{p}{k}}p
		-
		\logarithm{\sequence{p}{j}}\sequence{p}{k}
	}
	\\
	&
	=
	- \riemannian{X_{\sequence{p}{j}}}
	{
		(\cH_{t_0} - \id)[\logarithm{\sequence{p}{j}}\sequence{p}{k}]
	}
	\\
	&
	\le
	\riemanniannorm{X_{\sequence{p}{j}}}
	\norm{\cH_{t_0} - \id}
	\riemanniannorm{\logarithm{\sequence{p}{j}}\sequence{p}{k}}
	\\
	&
	\le
	\varrho_{\omega, \Omega}(\dist(\geodesic<s>(t_0), \, p))
	\riemanniannorm{X_{\sequence{p}{j}}}
	\riemanniannorm{\logarithm{\sequence{p}{j}}\sequence{p}{k}}
	,
\end{align*}
where the first equality comes is due to \zcref{eq:inner-product-hessian}, the first inequality comes from applying the Cauchy-Schwarz inequality, and the last inequality comes from \zcref{ineq:minmax-eigs}.
By rearranging the terms we get
\begin{multline}
	\label{eq:difference-of-logs-inequality}
	\riemannian{
		X_{\sequence{p}{j}}
	}
	{
		\logarithm{\sequence{p}{j}}p
		-
		\parallelTransport{\sequence{p}{k}}{\sequence{p}{j}}
		\logarithm{\sequence{p}{k}}p
	}
	\ge
	\riemannian{
		X_{\sequence{p}{j}}
	}
	{
		\logarithm{\sequence{p}{j}} \sequence{p}{k}
	}
	\\
	-
	\varrho_{\omega, \Omega}(\dist(\geodesic<s>(t_0), \, p))
	\riemanniannorm{X_{\sequence{p}{j}}}
	\riemanniannorm{\logarithm{\sequence{p}{j}}\sequence{p}{k}}
\end{multline}
for all $p \in \dom f$.
Since $\varrho_{\omega, \Omega}(\dist(\geodesic<s>(t_0), \, p)) \le \varrho_{\omega, \Omega}(\delta)$ for all $p \in \dom f$, we get
\begin{align*}
	\riemannian{
		X_{\sequence{p}{j}}
	}
	{
		\logarithm{\sequence{p}{j}}p
		-
		\parallelTransport{\sequence{p}{k}}{\sequence{p}{j}}
		\logarithm{\sequence{p}{k}}p
	}
	&
	\ge
	\riemannian{
		X_{\sequence{p}{j}}
	}
	{
		\logarithm{\sequence{p}{j}} \sequence{p}{k}
	}
	\\
	&
	\quad
	-
	\varrho_{\omega, \Omega}(\delta)
	\riemanniannorm{X_{\sequence{p}{j}}}
	\riemanniannorm{\logarithm{\sequence{p}{j}}\sequence{p}{k}}
	,
\end{align*}
for all $p \in \dom f$.
This is gives \zcref{ineq:almost-subgradient}.